\documentclass[11pt]{article}

\usepackage{amssymb}
\usepackage{amsmath}
\usepackage{amsthm}
\usepackage{mathrsfs}
\usepackage{graphicx}
\usepackage[utf8]{inputenc}
\usepackage{multirow}
\usepackage{floatrow}
\usepackage{mathrsfs}

\newtheorem{thm}{Theorem}[subsection]
\newtheorem{prop}[thm]{Proposition}
\newtheorem{lem}[thm]{Lemma}
\newtheorem{cor}[thm]{Corollary}

\newtheorem{defn}[thm]{Definition}

\newenvironment{Pf}{\noindent{\bf Proof. }}{\hfill $\blacksquare$ \\}

\allowdisplaybreaks

\def\lf{\lfloor}
\def\rf{\rfloor}
\def\lc{\lceil}
\def\rc{\rceil}
\def\ov{\overline}

\def\la{\langle}
\def\ra{\rangle}

\def\a{\alpha}
\def\b{\beta}
\def\e{\textbf e}

\def\g{\textbf g}
\def\k{\kappa}
\def\l{\lambda}
\def\m{\textbf m}
\def\r{\rho}
\def\s{\sigma}

\def\Ap{\text{Ap}}
\def\F{\textbf F}
\def\G{\text G}
\def\N{\mathbb N}

\def\R{\mathbb R}
\def\Z{\mathbb Z}

\def\ds{\displaystyle}
\def\sm{\setminus}
\def\pr{\prime}

\begin{document}

\title{\bf On the Frobenius Problem for Some Generalized Fibonacci Subsequences - II}
\author{
{\bf Ryan Azim Shaikh}\thanks{Department of Mathematics, Indian Institute of Technology, Hauz Khas, New Delhi -- 110016, India. \newline {\tt e-mail:mt1221923@maths.iitd.ac.in}} \qquad
{\bf Amitabha Tripathi}\thanks{Department of Mathematics, Indian Institute of Technology, Hauz Khas, New Delhi -- 110016, India. \newline {\tt e-mail:atripath@maths.iitd.ac.in}}\:\:\thanks{\it Corresponding author}
}
\date{}
\maketitle

\begin{abstract}
\noindent For a set $A$ of positive integers with $\gcd(A)=1$, let $\la A \ra$ denote the set of all finite linear combinations of elements of $A$ over the non-negative integers. Then it is well known that only finitely many positive integers do not belong to $\la A \ra$. The Frobenius number and the genus associated with the set $A$ is the largest number and the cardinality of the set of integers non-representable by $A$. By a generalized Fibonacci sequence $\{V_n\}_{n \ge 1}$ we mean any sequence of positive integers satisfying the recurrence $V_n=V_{n-1}+V_{n-2}$ for $n \ge 3$. We study the problem of determining  the Frobenius number and genus for sets $A=\{V_n, V_{n+d}, V_{n+2d}, \ldots \}$ for arbitrary $n$ and even $d$. 
\end{abstract}



\section{Introduction} \label{intro}
\vskip 10pt

For a given subset $A$ of positive integers with $\gcd(A)=1$, we write
\[ S = \la A \ra = \big\{ a_1x_1+\cdots+a_kx_k: a_i \in A, x_i \in {\Z}_{\ge 0} \big\}. \]
We say that $A$ is a set of generators for the set $S$. Further, $A$ is a minimal set of generators for $S$ if no proper subset of $A$ generates $S$. If $A=\{a_1, \ldots, a_n\}$ is any set of generators of $S$ arranged in increasing order, then $A$ is a minimal set of generators for $S$ if and only if $a_{k+1} \notin \la a_1, \ldots, a_k \ra$ for $k \in \{1,\ldots,n-1\}$. It is known that $A=S^{\star} \sm \big(S^{\star}+S^{\star}\big)$, where $S^{\star}=S \sm \{0\}$, is the unique minimal set of generators for $S$. The embedding dimension $\e(S)$ of $S$ is the size of the minimal set of generators. 

For any set of positive integers $A$ with $\gcd(A)=1$, the set ${\Z}_{\ge 0} \sm S$ is necessarily finite; we denote this by $\G(S)$. The cardinality of $\G(S)$ is the genus of $S$ and is denoted by $\g(S)$. The largest element in $\G(S)$ is the Frobenius number of $S$ and is denoted by $\F(S)$.

The Ap\'{e}ry set of $S$ corresponding to any fixed $a \in S$, denoted by $\Ap(S,a)$, consists of those $n \in S$ for which $n-a \notin S$. Thus, $\text{Ap}(S,a)$ is the set of minimum integers in $S \cap \textbf {C}$ as $\textbf{C}$ runs through the complete set of residue classes modulo $a$. 

The integers $\F(S)$ and $\g(S)$ can be computed from the Ap\'{e}ry set $\Ap(S,a)$ of $S$ corresponding to any $a \in S$ via the following proposition. 

\addtocounter{subsection}{1}

\begin{prop} {\bf (\cite{BS62,Sel77})} \label{prelims}
Let $S$ be a numerical semigroup, let $a \in S$, and let $\Ap(S,a)$ be the Ap\'{e}ry set of $S$ corresponding to $a$. Then 
\begin{itemize}
\item[{\rm (i)}]
\[ \F(S) = \max \Big( \Ap(S,a) \Big) - a; \]
\item[{\rm (ii)}]
\[ \g(S) = \frac{1}{a} \left( \sum_{n \in \Ap(S,a)} n \right) - \frac{a-1}{2}; \]
\end{itemize}  
\end{prop}
\vskip 5pt

The case where $\e(S)=2$ is well known and easy to establish. If $S=\la a,b \ra$, then it is easy to see that $\Ap(S,a)=\{bx: 0 \le x \le a-1\}$, and consequently 
\begin{equation} \label{e(S)=2}
\F(S) = ab-a-b, \quad \g(S) = \tfrac{1}{2} (a-1)(b-1) 
\end{equation} 
by Proposition \ref{prelims}. 
\vskip 5pt

The Frobenius Problem is the problem of determining the Frobenius number and the genus of a given numerical semigroup, and was first studied by Sylvester, and later by Frobenius; see \cite{Ram05} for a survey of the problem. Connections with Algebraic Geometry revived interest in Numerical Semigroups around the middle of the twentieth century; we refer to \cite{RG-S09} as a basic textbook on the subject. Curtis \cite{Cur90} proved that there exists no closed form expression for the Frobenius number of a numerical semigroup $S$ with $\e(S)>2$. As a consequence, a lot of research has focussed on the Frobenius number of semigroups whose generators are of a particular form. There are three particular instances of such results that are perhaps the closest to our work, and hence bear mentioning. Mar\'{i}n et. al. \cite{MRR07} determined the Frobenius number and genus of numerical semigroups of the form $\la F_i, F_{i+2}, F_{i+k} \ra$, where $i, k \ge 3$. These are called Fibonacci semigroups by the authors. Matthews \cite{Mat09} considers semigroups of the form $\la a, a+b, aF_{k-1}+bF_k \ra$ where $a>F_k$ and $\gcd(a,b)=1$. Taking $a=F_i$ and $b=F_{i+1}$, one gets the semigroup $\la F_i, F_{i+2}, F_{i+k} \ra$, considered in \cite{MRR07}. Thus, such semigroups were termed generalized Fibonacci semigroups by Matthews, who determined the Frobenius number of a generalized Fibonacci semigroup, thereby generalizing the result in \cite{MRR07} for Frobenius number. Batra et. al. \cite{BKT15} determined the Frobenius number and genus of numerical semigroups of the form $\la a, a+b, 2a+3b, \ldots, F_{2k-1}a+F_{2k}b \ra$ and $\la a, a+3b, 4a+7b, \ldots, L_{2k-1}a+L_{2k}b \ra$ where $\gcd(a,b)=1$. 

By a generalized Fibonacci sequence we mean any sequence $\{V_n\}$ of positive integers which satisfies the recurrence $V_n=V_{n-1}+V_{n-2}$ for each $n \ge 3$. A study of some subsequences of a generalized Fibonacci sequence $\{V_n\}$ was initiated by Panda et. al. \cite{PRT24}, in which the authors study the semigroup $S$ generated by $\la V_n, V_{n+d}, V_{n+2d}, \ra$ when $d$ is odd and when $d=2$, and $n$ is arbitrary. They show that $S$ is a numerical semigroup if and only if $\gcd(V_1,V_2)=1$ and $\gcd(V_n,F_d)=1$. The case of odd $d$ is easy to resolve since $\e(S)=2$, that is, each $V_{n+kd} \in \la V_n, V_{n+d} \ra$. For $d=2$, $\e(S)=\k$ where $\k$ satisfies $F_{2(\k-1)} \le V_n-1 < F_{2\k}$. Elements of the Ap\'{e}ry set $\Ap(S,V_n)$ are obtained by applying the Greedy Algorithm to each integer in $\{1,\ldots,V_n-1\}$ with respect to the sequence $F_2, F_4, F_6, \ldots $. There can be no closed form expression for this in general, but there is a simple expression for the Frobenius number in special cases $V_n=F_n$ and $V_n=L_n$, and a recurrence relation satisfied by $\g(S)$ in the special case $V_n=F_n$.  

This paper completes the study of the cases initiated in \cite{PRT24} by extending the results of $d=2$ to even $d$. Throughout this paper, let $S=\la V_n,V_{n+d},V_{n+2d}, \ldots \ra$, where $\gcd(V_1,V_2)=1$, $\gcd(V_n,F_d)=1$ and $d$ is even. The main results are similar to the ones in \cite{PRT24}; we list them below: 
\begin{itemize}
\item[{\rm (i)}] 
The embedding dimension $e(S)=\k$, where $\k$ is the smallest positive integer for which $F_{\k d}/F_d \ge V_n$; refer Theorem \ref{e(S)=k}. 
\item[{\rm (ii)}] 
The Ap\'{e}ry set $\Ap(S,V_n)=\left\{ V_{n+d}\,x - \left\lf \frac{F_{(k-1)d}\,x}{F_{kd}} \right\rf V_n: 1 \le x \le V_n-1 \right\} \cup \{0\}$; refer Theorem \ref{Apery_set} and Proposition \ref{d-representation}, part (iii). 
\item[{\rm (iii)}] 
The Frobenius number $\F(S)$ in the general case (refer Theorem \ref{Frob_gen}), and in the special cases when $V_n=F_n$ and $V_n=L_n$ (refer Corollary \ref{Frob_special}).  
\item[{\rm (iv)}]  
A recurrence for the genus $\g(S)$ in some special cases when $V_n=F_n$ and $V_n=L_n$ (refer Proposition \ref{genus_recurrence}).  
\end{itemize}
\vskip 10pt

\section{Preliminary Results} \label{prelim}
\vskip 10pt

A generalized Fibonacci sequence $\la V_n \ra_{n \ge 1}$ is defined by 
\begin{equation} \label{gen_def}
V_n = V_{n-1}+ V_{n-2}, \;\; n \ge 3, \quad \;\text{with}\; V_1 = a, V_2 = b, 
\end{equation}
where $a$ and $b$ are any positive integers. Two important special cases are (i) Fibonacci sequence $\{F_n\}_{n \ge 1}$ when $a=b=1$, and (ii) Lucas sequence $\{L_n\}_{n \ge 1}$ when $a=1$ and $b=3$. It is customary to extend these definitions to $F_0=F_2-F_1=0$ and $L_0=L_2-L_1=2$. Binet's formula give explicit values for $F_n$ and $L_n$: 
\[ F_n = \frac{{\a}^n-{\b}^n}{\a-\b}, \qquad L_n = {\a}^n+{\b}^n, \]
where $\a=(1+\sqrt{5})/2$ and $\b=(1-\sqrt{5})/2$ are the roots of the equation $x^2-x-1=0$. From these formulae, it is easy to see that $F_{2n}=L_nF_n$, and easy to derive 
\begin{equation} \label{FL_identities}
F_n^2 - F_{n+1} F_{n-1} = (-1)^{n-1}, \quad L_n^2 - L_{n+1} L_{n-1} = (-1)^n \cdot 5, \;\; n \ge 1. 
\end{equation}  
\vskip 5pt

\noindent The following bounds for $F_m V_n$ for the cases $V_n=F_n$ and $V_n=L_n$ when $d$ is even are useful in determining $\e(S)$; see Theorem \ref{e(S)=k}.  

\addtocounter{subsection}{1}

\begin{lem} \label{prod_ineq}
Let $m$ and $n$ be positive integers, with $m \ge n$.  
\begin{itemize}
\item[{\rm (i)}]
If $n \ge 3$, then 
\[ F_{m+n-2} < F_m F_n < F_{m+n-1}. \]
If $n=2$, then $F_{m+n-2}=F_m F_n<F_{m+n-1}$. \\
If $n=1$, then $F_{m+n-2}=F_m F_n$ if and only if $m=2$ and $F_m F_n=F_{m+n-1}$ holds for each $m$. 
\item[{\rm (ii)}]
If $n \ge 3$, then 
\[ F_{m+n-1} < F_m L_n < F_{m+n}. \]
If $n=2$, then $F_{m+n} \le F_m L_n < F_{m+n+1}$ and $F_m L_n=F_{m+n}$ if and only if $m=2$. \\
If $n=1$, then $F_m L_n=F_{m+n-1}$ holds for each $m$ and $F_m L_n=F_{m+n}$ if and only if $m=1$. 
\end{itemize}
\end{lem}

\begin{Pf}
\begin{itemize}
\item[{\rm (i)}]
The cases $n=1,2,3$ are easily verified. Assume the inequality holds for each positive integer $<n$, so that we have $F_{m+k-2}<F_m F_k<F_{m+k-1}$ for $k=n-1$ and $k=n-2$. Adding the two inequalities gives the desired inequality for $F_m F_n$. 
\item[{\rm (ii)}]
The cases $n=1,2,3$ are easily verified. Assume the inequality holds for each positive integer $<n$, so that we have $F_{m+k-1}<F_m L_k<F_{m+k}$ for $k=n-1$ and $k=n-2$. Adding the two inequalities gives the desired inequality for $F_m L_n$. 
\end{itemize}
\end{Pf}
\vskip 5pt
 
\noindent The following identities connecting generalized Fibonacci sequences with the Fibonacci sequence are useful in our subsequent work.  
\vskip 5pt

\begin{prop} \label{F_identity}
\begin{itemize}
\item[]



\item[{\rm (i)}]
For positive integers $m$ and $n$, 
\[ V_{m+n} = F_{n-1} V_m + F_n V_{m+1}. \]
In particular, $F_n \mid F_{kn}$ for each $k \ge 1$. 

\item[{\rm (ii)}]
For positive integers $m,n,d$, 
\[ F_n V_{m+n+d} - F_{n+d} V_{m+n} = (-1)^{n-1} F_d V_m. \]
In particular, for $k \ge1$, 
\[ F_d V_{n+kd} - F_{kd} V_{n+d} = (-1)^{d-1} F_{(k-1)d} V_n. \]

\item[{\rm (iii)}]
If $k \ge 2$, then 
\[ F_{kd} -  (L_d-1) F_{(k-1)d} = F_{(k-1)d} - F_{(k-2)d}. \]

\item[{\rm (iv)}]
If $k \ge 3$, then 
\[ \left\lf \frac{F_{kd}}{F_{(k-1)d}} \right\rf = L_d-1. \]

\item[{\rm (v)}]
If $k-1 \ge t \ge 2$, then 
\[ F_{kd} = (L_d-1) F_{(k-1)d} + (L_d-2) \sum_{i=t}^{k-2} F_{id} + (L_d-1) F_{(t-1)d} - F_{(t-2)d}. \]
In particular, 
\[ F_{kd} = (L_d-1) F_{(k-1)d} + (L_d-2) \sum_{i=2}^{k-2} F_{id} + (L_d-1) F_d. \]

\item[{\rm (vi)}]
If $k \ge t \ge 1$, then 
\[ (L_d-2) \sum_{i=t}^k V_{n+id} = \left( V_{n+(k+1)d} - V_{n+kd} \right) - \left( V_{n+td} - V_{n+(t-1)d} \right). \]
In particular, 
\[ (L_d-2) \sum_{i=1}^k V_{n+id} = \left( V_{n+(k+1)d} - V_{n+kd} \right) - \left( V_{n+d} - V_n \right). \]
\end{itemize}
\end{prop}
 
\begin{Pf}
\begin{itemize}


\item[{\rm (i)}]
We fix $m$ and induct on $n$. The case $n=1$ is an identity and the case $n=2$ follows from the definition of $\{V_n\}$. Assuming the result for all positive integers less than $n$, we have 
\begin{eqnarray*} 
V_{m+n} & = & V_{m+(n-1)} + V_{m+(n-2)} \\
& = & \big(F_{n-2} V_m + F_{n-1} V_{m+1}\big) +  \big(F_{n-3} V_m + F_{n-2} V_{m+1}\big) \\
& = & \big(F_{n-2}+F_{n-3}\big)V_m + \big(F_{n-1}+F_{n-2}\big)V_{m+1} \\
& = & F_{n-1} V_m + F_n V_{m+1}. 
\end{eqnarray*}
This completes the proof by induction. 

In particular, with $V_n=F_n$ and $m=(k-1)n$, we have  
\[ F_{kn} = F_{(k-1)n}F_{n-1} + F_{(k-1)n+1}F_n. \]
So if $F_n \mid F_{(k-1)n}$, then $F_n \mid F_{kn}$. Hence, $F_n \mid F_{kn}$ for each $k \ge 1$ by induction. 

\item[{\rm (ii)}]
We first prove the case $d=1$, then use this to prove the general case. 

By part (i) and eqn.~\eqref{FL_identities}, we have 
\begin{eqnarray*} 
F_n V_{m+n+1} - F_{n+1} V_{m+n} & = & F_n \big( F_n V_m + F_{n+1} V_{m+1} \big) - F_{n+1} \big( F_{n-1} V_m + F_n V_{m+1} \big) \\
& = & \big(F_n^2 - F_{n+1} F_{n-1} \big) V_m \\
& = & (-1)^{n-1} V_m. 
\end{eqnarray*}
This proves the case $d=1$. 

To prove the general case, by part (i), we have 
\[ V_{m+n+d} = F_{d-1} V_{m+n} + F_d V_{m+n+1} \;\;\text{and}\;\; F_{n+d} = F_{d-1} F_n + F_d F_{n+1}. \]
Therefore
\begin{eqnarray*}
F_n V_{m+n+d} - F_{n+d} V_{m+n} & = & F_n \left( F_{d-1} V_{m+n} + F_d V_{m+n+1} \right) - \left( F_{d-1} F_n + F_d F_{n+1} \right)  V_{m+n} \\
& = & F_d \left( F_n V_{m+n+1} - F_{n+1} V_{m+n} \right) \\
& = &  (-1)^{n-1} F_d V_m.
\end{eqnarray*}
This proves the general case. 
\vskip 5pt

Note that the particular case holds for $k=1$. For $k>1$, the transformation $m \mapsto n$, $n \mapsto d$, $d \mapsto (k-1)d$ yields the desired identity. 

\item[{\rm (iii)}]
Applying part (ii) to $V=F$ and $n=d$, and replacing $k$ by $k-1$, we have 
\begin{equation} \label{F_kd_eqn}
F_{kd} = \frac{F_{(k-1)d}}{F_d} F_{2d} - \frac{F_{(k-2)d}}{F_d} F_d = L_d F_{(k-1)d} - F_{(k-2)d} = (L_d-1) F_{(k-1)d} + \big( F_{(k-1)d} - F_{(k-2)d} \big). 
\end{equation}

\item[{\rm (iv)}]
Since $0<F_{(k-1)d} - F_{(k-2)d}<F_{(k-1)d}$ for $k \ge 3$, this follows upon dividing both sides of eqn.~\eqref{F_kd_eqn} by $F_{(k-1)d}$. 

\item[{\rm (v)}]
Replacing $k$ by $i$ in the identity in part (iii), then summing from $i=t$ to $i=k$, we have 
\[ \sum_{i=t}^k \left( F_{id} -  (L_d-1) F_{(i-1)d} \right) = \sum_{i=t}^k \left( F_{(i-1)d} - F_{(i-2)d} \right). \]
Thus, 
\[ \sum_{i=t}^k F_{id} -  (L_d-1) \sum_{i=t}^k F_{(i-1)d} = F_{kd} - (L_d-2) \sum_{i=t}^{k-1} F_{id} - (L_d-1) F_{(t-1)d} = F_{(k-1)d} - F_{(t-2)d}, \]
which gives the desired result. 

\item[{\rm (vi)}]
From part (v), 
\[ (L_d-2) \sum_{i=t}^k F_{id} = \left( F_{(k+1)d} - F_{kd} \right) - \left( F_{td} - F_{(t-1)d} \right). \]  
Replacing $k$ by $i$ in the identity in part (ii), then summing from $i=t$ to $i=k$ and multiplying both sides by $L_d-2$, we have 
\begin{eqnarray*} 
(L_d-2) \sum_{i=t}^k V_{n+id} & = & (L_d-2) \sum_{i=t}^k \frac{F_{id}}{F_d} V_{n+d} - (L_d-2) \sum_{i=t}^k \frac{F_{(i-1)d}}{F_d} V_n \\
& = & \Big( \left( F_{(k+1)d} - F_{kd} \right) - \left( F_{td} - F_{(t-1)d} \right) \Big) \frac{V_{n+d}}{F_d} \\
& & - \Big( \left( F_{kd} - F_{(k-1)d} \right) - \left( F_{(t-1)d} - F_{(t-2)d} \right) \Big) \frac{V_n}{F_d} \\ 
& = & \left( \frac{F_{(k+1)d}}{F_d} V_{n+d} - \frac{F_{kd}}{F_d} V_n \right) - \left( \frac{F_{kd}}{F_d} V_{n+d} - \frac{F_{(k-1)d}}{F_d} V_n \right) \\
& & - \left( \frac{F_{td}}{F_d} V_{n+d} - \frac{F_{(t-1)d}}{F_d} V_n \right) + \left( \frac{F_{(t-1)d}}{F_d} V_{n+d} - \frac{F_{(t-2)d}}{F_d} V_n \right) \\
& = & \left( V_{n+(k+1)d} - V_{n+kd} \right) - \left( V_{n+td} - V_{n+(t-1)d} \right).
\end{eqnarray*}
\end{itemize}
\end{Pf}
\vskip 5pt

\section{The Case where $d$ is even} \label{d_even}
\vskip 10pt

The main results of this paper are contained in this Section. We begin by proving an explicit formula for the embedding dimension in Theorem \ref{e(S)=k} in Subsection \ref{embedding_dim_subsection}. We follow this by introducing the Greedy Algorithm in Subsection \ref{greedy_subsection}, and apply it to compute a specific Ap\'{e}ry set in Subsection \ref{apery_subsection}. Finally, we compute the Frobenius number and genus in Subsection \ref{spcases_subsection} by using the results of the previous Subsections.  
\vskip 5pt
 
\subsection{Embedding Dimension} \label{embedding_dim_subsection}
\vskip 5pt

\begin{thm} \label{e(S)=k}
Let $S=\la V_n,V_{n+d},V_{n+2d}, \ldots \ra$, where $d$ is even and $\gcd(V_1,V_2)=\gcd(V_n,F_d)=1$. Then embedding dimension of $S$ is given by 
\[ e(S) = \k, \] 
where $\k$ is the smallest positive integer for which $F_{\k d}/F_d \ge V_n$. 
\end{thm}

\begin{Pf}
We claim that $\{V_n, V_{n+d}, V_{n+2d}, \ldots, V_{n+(\k -1)d}\}$ is a minimal set of generators for $S$, where $\k$ is the smallest positive integer for which $F_{\k d}/F_d \ge V_n$. By the characterization of minimal set of generators in Section \ref{intro}, we must therefore show: 
\begin{itemize}
\item[(i)]
$V_{n+kd} \in \la V_n, V_{n+d}, V_{n+2d}, \ldots, V_{n+(\k -1)d} \ra$ for each $k \ge \k$, and 
\item[(ii)]
$V_{n+kd} \notin \la V_n, V_{n+d}, V_{n+2d}, \ldots, V_{n+(k-1)d} \ra$ for $1 \le k \le \k-1$. 
\end{itemize}
\vskip 5pt

\noindent Let $k \ge \k$. By Proposition \ref{F_identity}, part (ii), we can write 
\begin{eqnarray} \label{identity1}
V_{n+kd} & = & -\frac{F_{(k-1)d}}{F_d}\, V_n + \frac{F_{kd}}{F_d}\, V_{n+d} \\
& = & \left( \l V_{n+d} - \frac{F_{(k-1)d}}{F_d} \right)V_n + \left( \frac{F_{kd}}{F_d} - \l V_n \right)V_{n+d} \;\;\text{ for any } \l \in \N. \nonumber
\end{eqnarray}
Therefore, $V_{n+kd} \in \la V_n, V_{n+d} \ra$ if there exists $\l \in \N$ for which 
\[ \frac{F_{(k-1)d}}{F_d\,V_{n+d}} \le \l \le \frac{F_{kd}}{F_d\,V_n}. \] 

\noindent If $F_{(k-1)d}/F_d\,V_{n+d} \le 1$, then $\l=1$ works because of the definition of $\k$. If $F_{(k-1)d}/F_d\,V_{n+d}>1$, then 
\[ \l\,V_n < \left( \frac{F_{(k-1)d}}{F_d\,V_{n+d}} + 1\right) V_n < 2\,\frac{F_{(k-1)d}}{F_d\,V_{n+d}}\,V_n < 2\,\frac{F_{(k-1)d}}{F_d} < \frac{F_{kd}}{F_d} \]
for $\l = \left\lc F_{(k-1)d}/F_d\,V_{n+d} \right\rc$ where the last inequality holds because $2\,F_{(k-1)d}<F_{kd}$. This proves claim (i). 
\vskip 5pt

\noindent Let $1<k<\k$. To prove claim (ii), suppose 
\begin{eqnarray} \label{identity2}
V_{n+kd} & = & \sum_{i=0}^{k-1} a_i V_{n+id} \nonumber \\
& = & a_0 \,V_n + a_1 \,V_{n+d} + \sum_{i=2}^{k-1} a_i \left( -\frac{F_{(i-1)d}}{F_d}\, V_n + \frac{F_{id}}{F_d}\, V_{n+d} \right) \nonumber \\
& = & \left( a_0 - \sum_{i=2}^{k-1} a_i \frac{F_{(i-1)d}}{F_d} \right) V_n + \left( a_1 + \sum_{i=2}^{k-1} a_i \frac{F_{id}}{F_d} \right) V_{n+d}.
\end{eqnarray} 
with each $a_i \ge 0$.

\noindent Note that $\gcd(V_n,V_{n+d})=1$ since any common divisor of $V_n$ and $V_{n+d}$ must divide each of the terms $V_{n+kd}$ due to eqn.~\eqref{identity1}. Thus, from eqn.~\eqref{identity1} and eqn.~\eqref{identity2}, there exists $t \in \Z$ such that 
\begin{eqnarray} 
\frac{F_{kd}}{F_d} + t\,V_n = a_1 + \sum_{i=2}^{k-1} a_i \frac{F_{id}}{F_d} \label{identity3} \\
-\frac{F_{(k-1)d}}{F_d} - t\,V_{n+d} = a_0 - \sum_{i=2}^{k-1} a_i \frac{F_{(i-1)d}}{F_d}. \label{identity4} 
\end{eqnarray}      
\vskip 5pt

\noindent In eqn.~\eqref{identity3}, $t<0$ reduces the left-side to a negative quantity, whereas the right-side is non-negative. Thus, $t \ge 0$. We rewrite eqn.~\eqref{identity3} and eqn.~\eqref{identity4} in the form 
\begin{eqnarray} 
\frac{F_{kd}}{F_d} \left( \sum_{i=2}^{k-1} a_i \frac{F_{id}}{F_{kd}} - 1 \right) = t\,V_n - a_1 \label{identity5} \\
\frac{F_{(k-1)d}}{F_d} \left( \sum_{i=2}^{k-1} a_i \frac{F_{(i-1)d}}{F_{(k-1)d}} - 1 \right) = t\,V_{n+d} + a_0. \label{identity6} 
\end{eqnarray}      
\vskip 5pt

\noindent With $m=(k-i)d$ and $n=(i-1)d$, and choosing $V=F$ in Proposition \ref{F_identity}, part (ii) we get 
\begin{equation} \label{inc_ratio}
\frac{F_{(i-1)d}}{F_{(k-1)d}} < \frac{F_{id}}{F_{kd}} 
\end{equation} 
for $1<i<k$. Using eqn.~\eqref{identity5} and eqn.~\eqref{identity6} now leads to the impossibility 
\[ t\,V_{n+d} + a_0 = \frac{F_{(k-1)d}}{F_d} \left( \sum_{i=2}^{k-1} a_i \frac{F_{(i-1)d}}{F_{(k-1)d}} - 1 \right) < \frac{F_{kd}}{F_d} \left( \sum_{i=2}^{k-1} a_i \frac{F_{id}}{F_{kd}} - 1 \right) = t\,V_n - a_1 \]
since $t \ge 0$. This proves claim (ii). 
\end{Pf}

\begin{cor} \label{emb_dim_fiblucas}
\begin{itemize}
\item[]
\item[{\rm (i)}]
If $d$ is even and $\gcd(F_n,F_d)=1$, the embedding dimension of $S_1=\la F_n, F_{n+d}, F_{n+2d}, \ldots \ra$ is given by 
\[ \e(S_1) = \begin{cases} 
1+\left\lc\frac{n-2}{d}\right\rc & \mbox{ if } d=2 \mbox{ or } d>2, n \le 2, \\ 
1+\left\lc\frac{n-1}{d}\right\rc & \mbox{ if } d>2, n>2.  
\end{cases} 
\]
\item[{\rm (ii)}]
If $d$ is even and $\gcd(L_n,F_d)=1$, the embedding dimension of $S_2=\la L_n, L_{n+d}, L_{n+2d}, \ldots \ra$ is given by 
\[ \e(S_2) = \begin{cases} 
1 & \mbox{ if } n=1, \\ 
1+\left\lc\frac{n}{d}\right\rc & \mbox{ if } n>1.  
\end{cases} 
\]
\end{itemize}
\end{cor}

\begin{Pf}
This is a direct application of Lemma \ref{prod_ineq} and Theorem \ref{e(S)=k}. 
\begin{itemize}
\item[{\rm (i)}]
If $d=2$ or $d>2$, $n \le 2$, then $\k$ is the least positive integer satisfying $\k d \ge n+d-2$. Hence $\k=1+\left\lc\frac{n-2}{d}\right\rc$ for these cases. \\
If $d>2$ and $n>2$, then $\k$ is the least positive integer satisfying $\k d \ge n+d-1$. Hence $\k=1+\left\lc\frac{n-1}{d}\right\rc$ for these cases.  
\item[{\rm (ii)}]
If $n=1$, then $\k$ is the least positive integer satisfying $\k d \ge d$. Hence $\k=1$ in this case. \\
If $n>1$, then $\k$ is the least positive integer satisfying $\k d \ge n+d$. Hence $\k=1+\left\lc\frac{n}{d}\right\rc$ for these cases.  
\end{itemize}
\end{Pf}
\vskip 5pt

\subsection{Some Results Based on the Greedy Algorithm} \label{greedy_subsection}
\vskip 5pt

\begin{defn} {\bf (The Greedy Algorithm)} \label{greedy_algo} \\[2pt]
For positive integers $c_1,\ldots,c_n,C$ with $\gcd(c_1,\ldots,c_n) \mid C$, consider the equation 
\begin{equation} \label{greedyeq}
c_1 x_1 + \cdots + c_n x_n = C. 
\end{equation}
The {\tt greedy solution\/} is given by 
\[ x_k^{\star} = \begin{cases}
                           \left\lf \frac{C}{c_n} \right\rf & \mbox{ for } k=n; \\[8pt]
                           \left\lf \frac{C - \sum_{i=k+1}^n c_i x_i^{\star}}{c_k} \right\rf & \mbox{ for } k=n-1, n-2, \ldots, 1.
                         \end{cases}
\]
We then write {\sc Greedy}$(c_1,\ldots,c_n;C)=x_1^{\star},\ldots,x_n^{\star}$. 
\end{defn}
\vskip 5pt

\begin{defn} \label{ess}
Fix $x \in \{1,\ldots,V_n-1\}$, and let $k$ be such that $F_{kd}/F_d \le x < F_{(k+1)d}/F_d$. Let 
\[ {\l}_1,\ldots,{\l}_k = \text{\sc Greedy}(1,F_{2d}/F_d,F_{3d}/F_d,\ldots,F_{kd}/F_d;x). \]
Set 
\[ s(x) = \sum_{i=1}^k {\l}_i V_{n+id}. \]
\end{defn}
\vskip 5pt

\begin{prop} \label{greedy_algo_for_X}
Let $d$ be even, $x$ be a positive integer and $k>1$. Suppose 
\[ \text{\sc Greedy}(1,F_{2d}/F_d,F_{3d}/F_d,\ldots,F_{kd}/F_d;x) = {\l}_1, \ldots, {\l}_k. \] 
\begin{itemize}
\item[{\rm (i)}]
Then $0 \le {\l}_i \le L_d-1$ for $1 \le i < k$. 
\item[{\rm (ii)}]
If ${\l}_i={\l}_j=L_d-1$ for some $i<j<k$, then ${\l}_t<L_d-2$ for some $t$ satisfying $i<t<j$. Moreover, there does not exist $i<k-1$ such that ${\l}_i={\l}_{i+1}=L_d-1$. 
\end{itemize}
\end{prop}

\begin{Pf}
\begin{itemize}
\item[{\rm (i)}]
We have 
\begin{equation} \label{greedy1}
{\l}_k = \left\lf \frac{x}{F_{kd}/F_d} \right\rf, \quad {\l}_j = \left\lf \frac{x-\sum_{i=j+1}^k {\l}_i F_{id}/F_d}{F_{jd}/F_d} \right\rf, \:\: j=k-1,k-2,\ldots,1. 
\end{equation}
by Definition \ref{greedy_algo}. 
\vskip 5pt

\noindent By eqn.~\eqref{greedy1} and Proposition \ref{F_identity}, part (iv), we have 
\[ {\l}_j = \left\lf \frac{x-\sum_{i=j+1}^k {\l}_i F_{id}/F_d}{F_{jd}/F_d} \right\rf \le  \left\lf \frac{F_{(j+1)d}/F_d}{F_{jd}/F_d} \right\rf = L_d-1, \] 
for $2 \le j \le k-1$, and 
\[ {\l}_1 = \left\lf \frac{x-\sum_{i=2}^k {\l}_i F_{id}/F_d}{F_d/F_d} \right\rf = x-\sum_{i=2}^k {\l}_i \frac{F_{id}}{F_d} < \frac{F_{2d}}{F_d} = L_d. \]
This completes the proof of part (i). 

\item[{\rm (ii)}]
Suppose ${\l}_i={\l}_j=L_d-1$ for some $i<j$ and ${\l}_t \ge L_d-2$ for $i<t<j$. Then 
\begin{equation} \label{contradiction}
x - \sum_{t=j+1}^k {\l}_t \frac{F_{td}}{F_d} \ge (L_d-2) \sum_{t=i}^j \frac{F_{td}}{F_d} + \frac{F_{id}}{F_d} + \frac{F_{jd}}{F_d} = \frac{F_{(j+1)d}}{F_d} + \frac{F_{(i-1)d}}{F_d} \ge \frac{F_{(j+1)d}}{F_d} 
\end{equation} 
using Proposition \ref{F_identity}, part (v). This contradicts the definition of ${\l}_{j+1}$.  
\end{itemize}

\noindent If ${\l}_i={\l}_{i+1}=L_d-1$ for some $i<k-1$, the argument in eqn.~\eqref{contradiction} with $j=i+1$ again leads to the same contradiction. This proves part (ii). 
\end{Pf}
\vskip 5pt

\begin{prop} \label{greedy_algo_for_F_kd}
Let $d$ be even and $k>1$. Suppose 
\[ \text{\sc Greedy}(1,F_{2d}/F_d,F_{3d}/F_d,\ldots,F_{kd}/F_d;F_{(k+1)d}/F_d) = {\l}_1, \ldots, {\l}_k. \] 
\begin{itemize}
\item[{\rm (i)}] 
\[ {\l}_i = \begin{cases}
                 L_d-1 & \mbox{ for } i=1, k; \\
                 L_d-2 & \mbox{ for } 2 \le i \le k-1.
               \end{cases}
\]
\item[{\rm (ii)}] 
\[ \text{\sc Greedy}(1,F_{2d}/F_d,F_{3d}/F_d,\ldots,F_{kd}/F_d;F_{(k+1)d}/F_d-1) = {\l}_1-1, \ldots, {\l}_k. \]
\item[{\rm (iii)}] 
\[ s\left(F_{(k+1)d}/F_d-1\right) = V_{n+(k+1)d} - V_{n+d} + V_n. \]
\end{itemize}
\end{prop}

\begin{Pf}
\begin{itemize}
\item[{\rm (i)}]
Observe that ${\l}_k=L_d-1$ follows from Proposition \ref{F_identity}, part (iv) and Definition \ref{greedy_algo}. 

\noindent We now prove that ${\l}_i=L_d-2$ for $2 \le i \le k-1$ by induction. We have 
\[ {\l}_{k-1} = \left\lf \frac{F_{(k+1)d} - (L_d-1)\,F_{kd}}{F_{(k-1)d}} \right\rf =  \left\lf \frac{F_{kd} - F_{(k-1)d}}{F_{(k-1)d}} \right\rf = L_d-2 \] 
from Proposition \ref{F_identity}, parts (iv) and (v), except that the last equality gives $L_d-1$ when $k=2$. 
 
\noindent Assuming ${\l}_j=L_d-2$ for some $j \in \{i+1,\ldots,k-1\}$, we have 
\[ {\l}_i = \left\lf \frac{F_{(k+1)d} - (L_d-1)\,F_{kd} - (L_d-2) \sum_{j=i+1}^{k-1} F_{jd}}{F_{id}} \right\rf = \left\lf \frac{F_{(i+1)d} - F_{id}}{F_{id}} \right\rf  = L_d-2 \] 
from Proposition \ref{F_identity}, parts (iv) and (v). 

\noindent Finally, we have 
\[ {\l}_1 = \left\lf \frac{F_{(k+1)d} - (L_d-1)\,F_{kd} - (L_d-2) \sum_{j=2}^{k-1} F_{jd}}{F_d} \right\rf = \left\lf \frac{F_{2d} - F_d}{F_d} \right\rf  = L_d-1 \] 
from Proposition \ref{F_identity}, part (v). 

\item[{\rm (ii)}]
Write $\text{\sc Greedy}(1,F_{2d}/F_d,F_{3d}/F_d,\ldots,F_{kd}/F_d;F_{(k+1)d}/F_d-1) = {\l}_1^{\star}, \ldots, {\l}_k^{\star}$. Then ${\l}_k^{\star}={\l}_k$ because $F_{kd} \nmid F_{(k+1)d}$ for $k>1$. Moreover, the numerator when computing ${\l}_i$ is $F_{(i+1)d}-F_{id}$; this is not a multiple of $F_{id}$ for $i>1$. Hence ${\l}_i^{\star}={\l}_i$ for $2 \le i \le k-1$. It follows that ${\l}_1^{\star}={\l}_1-1$.  

\item[{\rm(iii)}]
We have 
\begin{eqnarray*} 
s\left(F_{(k+1)d}/F_d-1\right) & = & \sum_{i=1}^k {\l}_i V_{n+id} - V_{n+d} \\
& = & (L_d-2) \sum_{i=1}^k V_{n+id} + V_{n+kd} \\
& = & \left( V_{n+(k+1)d} - V_{n+kd} - V_{n+d} + V_n \right) + V_{n+kd} \\
& = & V_{n+(k+1)d} - V_{n+d} + V_n. 
\end{eqnarray*} 
\end{itemize}
\end{Pf}
\vskip 5pt

\noindent We are now in a position to determine the Ap\'{e}ry set for the case $d$ even. We show that the elements in this set are obtained by applying the Greedy Algorithm to an equation involving terms of the form $F_{kd}$. 
\vskip 5pt
 
\begin{prop} \label{d-representation}
Fix $x \in \{1,\ldots,V_n-1\}$, and let $k$ be such that $F_{kd}/F_d \le x < F_{(k+1)d}/F_d$. Let 
\[ {\l}_1,\ldots,{\l}_k = \text{\sc Greedy}(1,F_{2d}/F_d,F_{3d}/F_d,\ldots,F_{kd}/F_d;x). \]
\begin{itemize}
\item[{\rm (i)}]
$0 \le {\l}_i \le L_d-1$ for each $i$ and ${\l}_k \ge 1$. 
\item[{\rm (ii)}]
$s(x) = \sum_{i=1}^k {\l}_i V_{n+id}$ satisfies 
\[ V_{n+kd} \le s(x) < V_{n+(k+1)d}, \quad s(x) \equiv V_{n+d}\,x\!\!\!\!\pmod{V_n}. \]
\item[{\rm (iii)}]
\[ s(x) = V_{n+d}\,x - \left\lf \frac{F_{(k-1)d}\,x}{F_{kd}} \right\rf V_n. \]
\end{itemize}
\end{prop}

\begin{Pf}
\begin{itemize}
\item[{\rm (i)}]
We define the sequence ${\l}_k, {\l}_{k-1}, \ldots, {\l}_1$ by using the Greedy Algorithm on $x$ with respect to the sequence $1, F_{2d}/F_d, F_{3d}/F_d, \ldots, F_{kd}/F_d$: 
\begin{equation} \label{greedy2}
{\l}_k = \left\lf \frac{x}{F_{kd}/F_d} \right\rf, \quad {\l}_j = \left\lf \frac{x-\sum_{i=j+1}^k {\l}_i F_{id}/F_d}{F_{jd}/F_d} \right\rf, \:\: j=k-1,k-2,\ldots,1. 
\end{equation}
\vskip 5pt

\noindent By Proposition \ref{F_identity}, part (iv), we have 
\[ 1 \le {\l}_k = \left\lf \frac{x}{F_{kd}/F_d} \right\rf \le \left\lf \frac{F_{(k+1)d}/F_d}{F_{kd}/F_d} \right\rf \le L_d-1, \] 
for $2 \le j \le k-1$, 
\[ {\l}_j = \left\lf \frac{x-\sum_{i=j+1}^k {\l}_i F_{id}/F_d}{F_{jd}/F_d} \right\rf \le  \left\lf \frac{F_{(j+1)d}/F_d}{F_{jd}/F_d} \right\rf \le L_d-1, \] 
and 
\[ {\l}_1 = \left\lf \frac{x-\sum_{i=2}^k {\l}_i F_{id}/F_d}{F_d/F_d} \right\rf = x-\sum_{i=2}^k {\l}_i \frac{F_{id}}{F_d} < \frac{F_{2d}}{F_d} = L_d. \]
This completes the proof of part (i). 

\item[{\rm (ii)}]
Define $s(x)=\sum_{i=1}^k {\l}_i V_{n+id}$. By Proposition \ref{F_identity}, part (ii), 
\[ s(x) \equiv \sum_{i=1}^k {\l}_i \frac{F_{id}}{F_d}\,V_{n+d} = V_{n+d}\,x \!\!\!\!\pmod{V_n}. \] 
Since ${\l}_k \ge 1$ and ${\l}_i \ge 0$ for $1 \le i \le k-1$, we have $s(x) \ge V_{n+kd}$. To prove the upper bound for $s(x)$, we consider two cases: (I) ${\l}_k \le L_d-2$, and (II) ${\l}_k=L_d-1$. 
\vskip 5pt

\noindent {\sc Case} (I): If ${\l}_k \le L_d-2$, then
\begin{eqnarray*} 
s(x) & \le & (L_d-1) \sum_{i=1}^{k-1} V_{n+id} + (L_d-2) V_{n+kd} \\
& = & \frac{L_d-1}{L_d-2} \Big( \big( V_{n+kd} - V_{n+(k-1)d} \big) - \big( V_{n+d} - V_n \big) \Big) + (L_d-2) V_{n+kd} \\  
& < & L_d V_{n+kd} - V_{n+(k-1)d} \\
& = & \frac{F_{2d}}{F_d} V_{n+kd} - \frac{F_d}{F_d} V_{n+(k-1)d} \\
& = & V_{n+(k+1)d} 
\end{eqnarray*} 
using Proposition \ref{F_identity}, parts (ii) and (vi). 
\vskip 5pt

\noindent {\sc Case} (II): Suppose ${\l}_k=L_d-1$. We claim that one of the following cases must arise: (i) ${\l}_i=L_d-2$ for $i \in \{1,\ldots,k-1\}$; (ii) there exists $r \in \{1,\ldots,k-1\}$ such that ${\l}_r<L_d-2$ and ${\l}_i=L_d-2$ for $i \in \{r+1,\ldots,k-1\}$. 
\vskip 5pt

\noindent If neither of these cases is true, then there must exist $t \in \{1,\ldots,k-1\}$ such that ${\l}_t=L_d-1$ and ${\l}_i=L_d-2$ for $i \in \{t+1,\ldots,k-1\}$. But then 
\[ x \ge (L_d-2) \sum_{i=t}^k \frac{F_{id}}{F_d} + \frac{F_{td}}{F_d} + \frac{F_{kd}}{F_d} = \frac{F_{(k+1)d}}{F_d} + \frac{F_{(t-1)d}}{F_d} \ge \frac{F_{(k+1)d}}{F_d} \] 
using Proposition \ref{F_identity}, part (v). This contradiction proves the claim. 
\vskip 5pt

\noindent In case (i), we have 
\[ s(x) = (L_d-2) \sum_{i=1}^k V_{n+id} + V_{n+kd} = V_{n+(k+1)d} - \big( V_{n+d} - V_n \big) < V_{n+(k+1)d} \]
using Proposition \ref{F_identity}, part (vi). 

\noindent In case (ii), we have 
\begin{eqnarray*} 
s(x) & \le & (L_d-2) \sum_{i=r}^k V_{n+id} + (L_d-1) \sum_{i=1}^{r-1} V_{n+id} + V_{n+kd} - V_{n+rd} \\
& = & \big(V_{n+(k+1)d} - V_{n+kd}\big) - \big( V_{n+rd} - V_{n+(r-1)d} \big) \\ 
& & + \frac{L_d-1}{L_d-2} \left( \big(V_{n+rd} - V_{n+(r-1)d}\big) - \big( V_{n+d} - V_n \big) \right) + V_{n+kd} - V_{n+rd} \\
& = & V_{n+(k+1)d} - \frac{L_d-3}{L_d-2} V_{n+rd} - \frac{1}{L_d-2} V_{n+(r-1)d} - \frac{L_d-1}{L_d-2} \big( V_{n+d} - V_d \big) \\
& < & V_{n+(k+1)d} 
\end{eqnarray*}
using Proposition \ref{F_identity}, part (viii). This completes the proof of part (ii), and the case (II) of the Proposition. 

\item[{\rm (iii)}]
By part (ii), we know that $s(x)=V_{n+d}\,x - \l\,V_n$ for some integer $\l$; we must show that $\l=\left\lf \frac{F_{(k-1)d}\,x}{F_{kd}} \right\rf$.
\vskip 5pt

Applying Proposition \ref{F_identity}, part (ii) we have  
\begin{eqnarray*}
s(x) & = & \sum_{i=1}^k {\l}_i V_{n+id} \\
& = & \sum_{i=1}^k {\l}_i \left( \frac{F_{id}}{F_d} V_{n+d} - \frac{F_{(i-1)d}}{F_d} V_n \right) \\
& = & V_{n+d}\,x - \left( \sum_{i=1}^k {\l}_i \frac{F_{(i-1)d}}{F_d} \right) V_n. 
\end{eqnarray*}
Therefore, we must show that 
\begin{equation}  \label{ess_formula}
\sum_{i=1}^k {\l}_i \frac{F_{(i-1)d}}{F_d} = \left\lf \frac{F_{(k-1)d}\,x}{F_{kd}} \right\rf.  
\end{equation}
\vskip 5pt

Using eqn.~\eqref{inc_ratio}, we have 
\[ \sum_{i=1}^k {\l}_i \frac{F_{(i-1)d}}{F_d}  = \sum_{i=1}^k {\l}_i \frac{F_{(i-1)d}}{F_{id}} \cdot \frac{F_{id}}{F_d} \le \frac{F_{(k-1)d}}{F_{kd}} \sum_{i=1}^k {\l}_i \frac{F_{id}}{F_d} = \frac{F_{(k-1)d}\,x}{F_{kd}}. \]
\end{itemize}
\vskip 5pt

\noindent Thus, to prove eqn.~\eqref{ess_formula}, we must show 
\[  \sum_{i=1}^k {\l}_i \frac{F_{(i-1)d}}{F_d} > \frac{F_{(k-1)d}\,x}{F_{kd}} - 1 = \frac{F_{(k-1)d}}{F_{kd}} \sum_{i=1}^k {\l}_i \frac{F_{id}}{F_d} - 1, \] 
which is equivalent to 
\[ \sum_{i=1}^k {\l}_i \left( \frac{F_{(k-1)d}}{F_{kd}} \cdot \frac{F_{id}}{F_d} - \frac{F_{(i-1)d}}{F_d} \right) < 1, \] 
and hence to 
\[ \sum_{i=1}^k {\l}_i \frac{F_{(k-i)d}}{F_{kd}} < 1 \]
by Proposition \ref{F_identity}, part (ii). 
\vskip 5pt

\noindent To prove the above inequality, we consider two cases: (I) ${\l}_1 \le L_d-2$, and (II) ${\l}_1=L_d-1$. The argument is along the same lines as for the upper bound in part (ii). 
\vskip 5pt

\noindent {\sc Case} (I): If ${\l}_1 \le L_d-2$, then  
\begin{eqnarray*} 
\sum_{i=1}^k {\l}_i F_{(k-i)d} & \le & (L_d-1) \sum_{i=1}^{k-1} F_{id} - F_{(k-1)d} \\
& = & \frac{L_d-1}{L_d-2} \big( F_{kd} - F_{(k-1)d} - F_d \big) - F_{(k-1)d} \\  
& < & F_{kd} + \frac{F_{kd} - (2L_d-3) F_{(k-1)d}}{L_d-2} \\
& < &  F_{kd} - \frac{(L_d-3) F_{(k-1)d}}{L_d-2} \\
& < & F_{kd}.
\end{eqnarray*} 
using Proposition \ref{F_identity}, parts (ii) and (v). 
\vskip 5pt

\noindent {\sc Case} (II): Suppose ${\l}_1=L_d-1$. We claim that one of the following cases must arise: (i) ${\l}_i=L_d-2$ for $i \in \{2,\ldots,k-1\}$; (ii) there exists $r \in \{2,\ldots,k-1\}$ such that ${\l}_r<L_d-2$ and ${\l}_i=L_d-2$ for $i \in \{1,\ldots,r-1\}$. 
\vskip 5pt

\noindent If neither of these cases is true, then there must exist $t \in \{2,\ldots,k-1\}$ such that ${\l}_t=L_d-1$ and ${\l}_i=L_d-2$ for $i \in \{2,\ldots,t-1\}$. But then 
\[ {\l}_{t+1} = \left\lf \frac{x-\sum_{i=t+2}^k {\l}_i F_{id}/F_d}{F_{(t+1)d}/F_d} \right\rf  = \left\lf \frac{\sum_{i=1}^{t+1} {\l}_i F_{id}}{F_{(t+1)d}} \right\rf  = {\l}_{t+1} + \left\lf \frac{(L_d-2) \sum_{i=1}^t F_{id}+F_d+F_{td}}{F_{(t+1)d}} \right\rf = {\l}_{t+1}+1 \] 
using Proposition \ref{F_identity}, part (v). This contradiction proves the claim. 
\vskip 5pt

\noindent In case (i), we have 
\[ \sum_{i=1}^{k-1} {\l}_i F_{(k-i)d} = (L_d-2) \sum_{i=1}^{k-1} F_{id} + F_{(k-1)d} = F_{kd} - F_d < F_{kd} \]
using Proposition \ref{F_identity}, part (v). 

\noindent In case (ii), we have 
\begin{eqnarray*} 
\sum_{i=1}^{k-1} {\l}_i F_{(k-i)d} & \le & (L_d-2) \sum_{i=1}^r F_{(k-i)d} + (L_d-1) \sum_{i=r+1}^{k-1} F_{(k-i)d} + F_{(k-1)d} - F_{(k-r)d} \\
& = & \big( F_{kd} - F_{(k-1)d} - F_{(k-r)d} + F_{(k-1-r)d} \big) + \frac{L_d-1}{L_d-2} \big( F_{(k-r)d} - F_{(k-1-r)d} - F_d \big) \\
& & + F_{(k-1)d} - F_{(k-r)d} \\ 
& = & F_{kd} - \frac{L_d-3}{L_d-2} F_{(k-r)d} - \frac{1}{L_d-2} F_{(k-1-r)d} - \frac{L_d-1}{L_d-2} F_d  \\
& < & F_{kd} 
\end{eqnarray*}
using Proposition \ref{F_identity}, part (v). 
\vskip 5pt

This completes the proof of eqn.~\eqref{ess_formula}, and of part (iii). 
\end{Pf}
\vskip 5pt

\begin{prop} \label{greedy_fibonacci}
Let $d$ be even and let $F_{kd}/F_d \le F_m-1 < F_{(k+1)d}/F_d$, $m>2$. Suppose 
\[ \text{\sc Greedy}(1,F_{2d}/F_d,F_{3d}/F_d,\ldots,F_{kd}/F_d;F_m-1) = {\l}_1, \ldots, {\l}_k. \]
\begin{itemize}
\item[{\rm (i)}]
If $m \equiv r\pmod{d}$, $m$ is odd, $1 \le r \le d-1$, then 
\[ {\l}_i = 
\begin{cases} 
L_d-2 & \mbox{ if } 1 \le i \le \frac{m-r}{d}-1; \\
L_d-F_{d-r}-1 & \mbox{ if } i=\frac{m-r}{d}; \\
F_r-1 & \mbox{ if } i=\frac{m-r}{d}+1. 
\end{cases}
\]
\item[{\rm (ii)}]
If $m \equiv r\pmod{d}$, $m \ge d$, $m$ is even, $0 \le r \le d-1$, then 
\[ {\l}_i = 
\begin{cases} 
L_d-2 & \mbox{ if } 1 \le i \le \frac{m-r}{d}-2; \\
L_d-1 & \mbox{ if } i=\frac{m-r}{d}-1; \\
F_{d-r}-1 & \mbox{ if } i=\frac{m-r}{d}; \\
F_r & \mbox{ if } i=\frac{m-r}{d}+1. 
\end{cases}
\]
If $m$ is even and $m<d$, then $k=1$ and ${\l}_1=F_m-1$. 

\item[{\rm (iii)}]
\[ s(F_m-1) = F_d V_{n+m} - V_{n+d} + V_n. \]
\end{itemize}
\end{prop}

\begin{Pf}
Let $m=qd+r$, $0 \le r \le d-1$. If $r=0$, then $\gcd(V_m,F_d)=F_d$ for the case $V_i=F_i$. So in order for $\la V_n, V_{n+d}, V_{n+2d}, \ldots \ra$ to exist, we may consider only the case $d=2$. We break up the case $r>0$ according as $r$ is odd or even in parts (i) and (ii), and deal with the case $r=0$ in part (iii). 

From Proposition \ref{F_identity}, part (ii) replacing $m$ by $n+qd$, $n$ by $r$, and $d$ by $d-r$ gives 
\begin{equation} \label{V_identity}
F_d V_{n+qd+r} = F_r V_{n+(q+1)d} + (-1)^r F_{d-r} V_{n+qd}. 
\end{equation}

\begin{itemize}
\item[{\rm (i)}]
We first claim that $q=k$ if $r=1$ and $q=k-1$ if $r>1$. 

If $r=1$, then 
\[  \frac{F_{qd}}{F_d} \le \left( F_{d+1} -1 \right) \frac{F_{qd}}{F_d} + \frac{F_{qd}-F_{(q-1)d}}{F_d} - 1 = F_{d+1} \frac{F_{qd}}{F_d} - \frac{F_{(q-1)d}}{F_d} - 1 = F_m-1 < \frac{F_{(q+1)d}}{F_d} \]
by Proposition \ref{F_identity}, part (ii) and eqn.~\eqref{V_identity} with $n=0$ and $V=F$. Thus, $q=k$ and ${\l}_k=\left\lf \frac{F_m-1}{F_{qd}/F_d} \right \rf=F_{d+1}-1=L_d-F_{d-1}-1$. Observe that ${\l}_{q+1}$ does not exist, ${\l}_q=L_d-F_{d-r}-1$ and $F_r-1=0$ in this case. 

If $r>1$, then 
\[  \frac{F_{(q+1)d}}{F_d} \le \left( F_r -1 \right) \frac{F_{(q+1)d}}{F_d} + \frac{F_{(q+1)d}-F_{d-r}\,F_{qd}}{F_d} - 1 = F_r \frac{F_{(q+1)d}}{F_d} - F_{d-r} \frac{F_{qd}}{F_d} - 1 = F_m-1 \]
and 
\[ F_m-1 < F_r \frac{F_{(q+1)d}}{F_d}  < \frac{F_{(q+2)d}}{F_d} \]
by eqn.~\eqref{V_identity} with $n=0$ and $V=F$ since $\frac{F_{(j+1)d}}{F_{jd}} \ge L_d-1 > F_d$ by Proposition \ref{F_identity}, part (iv). Thus, $q+1=k$ and ${\l}_k=\left\lf \frac{F_m-1}{F_{(q+1)d}/F_d} \right \rf=F_r-1$. 
Furthermore, 
\begin{eqnarray*} 
F_m - 1 - {\l}_{q+1} \frac{F_{(q+1)d}}{F_d} & = & \frac{F_{(q+1)d}-F_{d-r}\,F_{qd}}{F_d} - 1 \\
& = & \left( L_d - F_{d-r} \right) \frac{F_{qd}}{F_d} - \frac{F_{(q-1)d}}{F_d} - 1 \\
& = & \left( L_d - F_{d-r} - 1 \right) \frac{F_{qd}}{F_d} + \frac{F_{qd}-F_{(q-1)d}}{F_d} - 1
\end{eqnarray*}  
by Proposition \ref{F_identity}, part (ii). Thus, ${\l}_q=L_d-F_{d-r}-1$. 

We have shown that ${\l}_q=L_d-F_{d-r}-1$ in both cases, and ${\l}_{q+1}=F_r-1$ for $r>1$. We now show that ${\l}_i=L_d-2$ for $1 \le i \le q-1$ by induction. Now 
\begin{eqnarray*}
{\l}_{q-1} & = & \left\lf \frac{(F_m-1)-{\l}_{q+1}F_{(q+1)d}/F_d-{\l}_q F_{qd}/F_d}{F_{(q-1)d}/F_d} \right\rf \\
& = & \left\lf \frac{(F_{qd}-F_{(q-1)d})/F_d - 1}{F_{(q-1)d}/F_d} \right\rf \\
& = & \left\lf \frac{F_{qd}-F_d}{F_{(q-1)d}} \right\rf - 1 \\
& = & \left\lf \frac{L_d F_{(q-1)d} - F_{(q-2)d} - F_d}{F_{(q-1)d}} \right\rf - 1 \\
& = & L_d-2. 
\end{eqnarray*}

Suppose ${\l}_i=L_d-2$ and $(F_m-1)-\sum_{j=i}^{q+1} {\l}_j F_{jd}/F_d = (F_{id}-F_{(i-1)d})/F_d -1$ for some $i \le q-1$. We must show that ${\l}_{i-1}=L_d-2$. We have
\begin{eqnarray*}
{\l}_{i-1} & = & \left\lf \frac{(F_m-1)-\sum_{j=i}^{q+1} {\l}_j F_{jd}/F_d}{F_{(i-1)d}/F_d} \right\rf \\
& = & \left\lf \frac{(F_{id}-F_{(i-1)d})/F_d - 1}{F_{(i-1)d}/F_d} \right\rf \\
& = & \left\lf \frac{F_{id}-F_d}{F_{(i-1)d}} \right\rf - 1 \\
& = & \left\lf \frac{L_d F_{(i-1)d} - F_{(i-2)d} - F_d}{F_{(i-1)d}} \right\rf - 1 \\
& = & L_d-2. 
\end{eqnarray*} 

\item[{\rm (ii)}]
If $m<d$, then $F_m-1<F_d<L_d=F_{2d}/F_d$. Thus, $k=1$ and ${\l}_1=F_m-1$. 

Let $m \ge d$, so $q \ge 1$. We consider two cases: (a) $r=0$, (b) $1 \le r \le d-1$. 

\noindent {\sc Case} (a). For $d=2$, since $F_2=1$ and $F_{2(q-1)}<F_{2q}-1<F_{2q}$, we have $q=k+1$ and ${\l}_{q-1}=\left\lf \frac{F_{2q}-1}{F_{2(q-1)}}\right\rf=2$. The proof that ${\l}_i=1$ for $1 \le i \le q-2$ follows along the same lines by induction as in case (i); note that $L_2-2=1$. 
\vskip 5pt

Let $d>2$. We have 
\[ \frac{F_{qd}}{F_d} \le F_{qd} - 1 = F_m - 1 < \frac{F_{(q+1)d}}{F_d} \]
since $\frac{F_{(q+1)d}}{F_{qd}} \ge L_d-1 > F_d$ by Proposition \ref{F_identity}, part (iv). Thus, $q=k$ and ${\l}_k=\left\lf \frac{F_m-1}{F_{qd}/F_d} \right \rf=F_d-1$. 

Furthermore, 
\[  F_m - 1 - {\l}_q \frac{F_{qd}}{F_d} = \frac{F_{qd}}{F_d} - 1. \] 
Therefore, ${\l}_{k-1}, \ldots, {\l}_1$ are determined by Proposition \ref{greedy_algo_for_F_kd}, parts (i), (ii). Observe that ${\l}_{q+1}$ does not exist since $F_r=0$ in this case. 
\vskip 5pt

\noindent {\sc Case} (b). In this case, we have  
\[ \frac{F_{(q+1)d}}{F_d} \le F_r \frac{F_{(q+1)d}}{F_d} + F_{d-r} \frac{F_{qd}}{F_d} - 1 = F_m-1 \]
and 
\[ F_m-1 = (F_r+1) \frac{F_{(q+1)d}}{F_d} - \frac{F_{(q+1)d}-F_{d-r}\,F_{qd}}{F_d} - 1 < (F_r+1) \frac{F_{(q+1)d}}{F_d} < \frac{F_{(q+2)d}}{F_d} \]
by eqn.~\eqref{V_identity} with $n=0$ and $V=F$ since $\frac{F_{(j+1)d}}{F_{jd}} \ge L_d-1 > F_d$ by Proposition \ref{F_identity}, part (iv). Thus, $q+1=k$ and ${\l}_k=\left\lf \frac{F_m-1}{F_{(q+1)d}/F_d} \right \rf=F_r$. 
Furthermore, 
\[  F_m - 1 - {\l}_{q+1} \frac{F_{(q+1)d}}{F_d} = \frac{F_{d-r}\,F_{qd}}{F_d} - 1 = \left( F_{d-r} - 1 \right) \frac{F_{qd}}{F_d} + \frac{F_{qd}}{F_d} - 1. \]
Thus, ${\l}_q=\left\lf \frac{(F_m-1)-{\l}_{q+1}F_{(q+1)d}/F_d}{F_{qd}/F_d}\right\rf=F_{d-r}-1$. Now 
\begin{eqnarray*} 
F_m - 1 - {\l}_{q+1} \frac{F_{(q+1)d}}{F_d} - {\l}_q \frac{F_{qd}}{F_d} & = & \frac{F_{qd}}{F_d} - 1 \\
& = & L_d \frac{F_{(q-1)d}}{F_d} - \frac{F_{(q-2)d}}{F_d} - 1 \\
& = & (L_d - 1) \frac{F_{(q-1)d}}{F_d} + \frac{F_{(q-1)d}-F_{(q-2)d}}{F_d} - 1
\end{eqnarray*} 
by Proposition \ref{F_identity}, part (ii). Thus, ${\l}_{q-1}=\left\lf \frac{(F_m-1)-{\l}_{q+1}F_{(q+1)d}/F_d-{\l}_q F_{qd}/F_d}{F_{(q-1)d}/F_d}\right\rf=L_d-1$.
\vskip 5pt

We have shown that ${\l}_{q+1}=F_r$, ${\l}_q=F_{d-r}-1$ and ${\l}_{q-1}=L_d-1$. The proof for ${\l}_i=L_d-2$, $1 \le i \le q-2$, by induction is identical to the one provided in case (i), for odd $m$. 

\item[{\rm (iii)}]
In case (i), using eqn.~\eqref{V_identity}, we have  
\begin{eqnarray*}
s(F_m-1) & = & \sum_{i=1}^k {\l}_i V_{n+id} \\
& = & (L_d-2) \sum_{i=1}^{q-1} V_{n+id} + (L_d-F_{d-r}-1) V_{n+qd} + (F_r-1) V_{n+(q+1)d} \\
& = & \left( V_{n+qd} - V_{n+(q-1)d} - V_{n+d} + V_n \right) + (L_d-F_{d-r}-1) V_{n+qd} + (F_r-1) V_{n+(q+1)d} \\
& = & \left( L_d V_{n+qd} - V_{n+(q-1)d} \right) - V_{n+d} + V_n - F_{d-r} V_{n+qd} +  (F_r-1) V_{n+(q+1)d} \\
& = & V_{n+(q+1)d} - V_{n+d} + V_n - F_{d-r} V_{n+qd} +  (F_r-1) V_{n+(q+1)d} \\
& = & \left( F_r V_{n+(q+1)d} - F_{d-r} V_{n+qd} \right) - V_{n+d} + V_n \\
& = & F_d V_{n+qd+r} - V_{n+d} + V_n.  
\end{eqnarray*}
In case (ii), using eqn.~\eqref{V_identity}, we have  
\begin{eqnarray*}
s(F_m-1) & = & \sum_{i=1}^k {\l}_i V_{n+id} \\
& = & (L_d-2) \sum_{i=1}^{q-1} V_{n+id} + V_{n+(q-1)d} + (F_{d-r}-1) V_{n+qd} + F_r V_{n+(q+1)d} \\
& = & \left( V_{n+qd} - V_{n+(q-1)d} - V_{n+d} + V_n \right) + V_{n+(q-1)d} + (F_{d-r}-1) V_{n+qd} + F_r V_{n+(q+1)d} \\
& = & \left( F_r V_{n+(q+1)d} + F_{d-r} V_{n+qd} \right) - V_{n+d} + V_n \\ 
& = & F_d V_{n+qd+r} - V_{n+d} + V_n.  
\end{eqnarray*}
In case (iii), 
\begin{eqnarray*}
s(F_m-1) & = & \sum_{i=1}^k {\l}_i V_{n+2i} \\
& = & \sum_{i=1}^{q-1} V_{n+2i} + V_{n+2(q-1)} \\
& = & \left( V_{n+2q} - V_{n+2(q-1)} - V_{n+2} + V_n \right) + V_{n+2(q-1)} \\
& = & V_{n+2q} - V_{n+2} + V_n.  
\end{eqnarray*}
In each case, $s(F_m-1)=F_d V_{n+m}-V_{n+d}+V_n$. 
\end{itemize}
\end{Pf}
\vskip 5pt

\begin{prop} \label{greedy_lucas}
Let $d$ be even and let $F_{kd}/F_d \le L_m-1 < F_{(k+1)d}/F_d$, $m>1$. Suppose 
\[ \text{\sc Greedy}(1,F_{2d}/F_d,F_{3d}/F_d,\ldots,F_{kd}/F_d;L_m-1) = {\l}_1, \ldots, {\l}_k. \]
\begin{itemize}
\item[{\rm (i)}]
If $m \equiv r\pmod{d}$, $m$ is odd, $1 \le r \le d-1$, then 
\[ {\l}_i = 
\begin{cases} 
L_d-2 & \mbox{ if } 1 \le i \le \frac{m-r}{d}-2; \\
L_d-1 & \mbox{ if } i=\frac{m-r}{d}-1; \\
L_{d-r}-1 & \mbox{ if } i=\frac{m-r}{d}; \\
L_r & \mbox{ if } i=\frac{m-r}{d}+1. 
\end{cases}
\]
\item[{\rm (ii)}]
If $m \equiv r\pmod{d}$, $m$ is even, $1 \le r \le d-1$, then 
\[ {\l}_i = 
\begin{cases} 
L_d-2 & \mbox{ if } 1 \le i \le \frac{m-r}{d}-1; \\
L_d-L_{d-r}-1 & \mbox{ if } i=\frac{m-r}{d}; \\
L_r-1 & \mbox{ if } i=\frac{m-r}{d}+1. 
\end{cases}
\]
\item[{\rm (iii)}]
If $d \mid m$, then 
\[ {\l}_i = 
\begin{cases} 
L_d-2 & \mbox{ if } 1 \le i \le \frac{m}{d}-2; \\
L_d-3 & \mbox{ if } i=\frac{m}{d}-1; \\
L_d-1 & \mbox{ if } i=\frac{m}{d}. 
\end{cases}
\]
\item[{\rm (iv)}]
\[ s(L_m-1) = F_d \left(V_{n+m+1} + V_{n+m-1} \right) - V_{n+d} + V_n. \]
\end{itemize}
\end{prop}

\begin{Pf}
Let $m=dq+r$, $0 \le r \le d-1$. We break up the case $r>0$ according as $r$ is odd or even in parts (i) and (ii), and deal with the case $r=0$ in part (iii). 
\vskip 5pt

Recall $F_{r+1}+F_{r-1}=L_r$ for $r \ge 1$. Replacing $r$ in eqn.~\eqref{V_identity} first by $r+1$ and then by $r-1$, and adding the two resultant equations, with $V_i=F_i$ and $n=0$, we have     
\begin{equation} \label{L_identity}
F_d L_{qd+r} = L_r F_{(q+1)d} - (-1)^r L_{d-r} F_{qd}. 
\end{equation}

The proof of parts (i), (ii) follow on exactly the same lines as the proof of Proposition \ref{greedy_fibonacci}, the esssential difference being in the use of the identity in eqn.~\eqref{L_identity}.
\vskip 5pt

\noindent We prove part (iii). We have  
\[ \frac{F_{qd}}{F_d} \le (L_d-1) \frac{F_{qd}}{F_d} + \frac{F_{qd}-2\,F_{(q-1)d}}{F_d} - 1 = L_d \frac{F_{qd}}{F_d} - L_0 \frac{F_{(q-1)d}}{F_d} - 1 = L_{qd} - 1 \]
and 
\[ L_{qd} - 1 = L_d \frac{F_{qd}}{F_d} - L_0 \frac{F_{(q-1)d}}{F_d} - 1 < \frac{L_d\,F_{qd}-F_{(q-1)d}}{F_d} = \frac{F_{(q+1)d}}{F_d} \]
by Proposition \ref{F_identity}, part (ii). Thus, $q=k$ and ${\l}_q=\left\lf \frac{L_m-1}{F_{qd}/F_d} \right \rf=L_d-1$. Furthermore, 
\[ {\l}_{q-1} = \left\lf \frac{(L_m-1)-{\l}_q\,F_{qd}/F_d}{F_{(q-1)d}/F_d}\right\rf = \left\lf \frac{F_{qd}-2\,F_{(q-1)d}-F_d}{F_{(q-1)d}} \right\rf = L_d-3 \] 
by Proposition \ref{F_identity}, part (iv). 

We now show that ${\l}_i=L_d-2$ for $1 \le i \le q-2$ by induction. Now 
\[ {\l}_{q-2} = \left\lf \frac{(L_m-1)-{\l}_q\,F_{qd}/F_d-{\l}_{q-1} F_{(q-1)d}/F_d}{F_{(q-2)d}/F_d} \right\rf = \left\lf \frac{F_{(q-1)d}-F_{(q-2)d}-F_d}{F_{(q-2)d}} \right\rf = L_d-2  \] 
by Proposition \ref{F_identity}, part (iv). 

Suppose ${\l}_i=L_d-2$ and $(L_m-1)-\sum_{j=i}^q {\l}_j F_{jd}/F_d = (F_{id}-F_{(i-1)d})/F_d -1$ for some $i \le q-2$. We must show that ${\l}_{i-1}=L_d-2$. We have
\begin{eqnarray*}
{\l}_{i-1} & = & \left\lf \frac{(L_m-1)-\sum_{j=i}^q {\l}_j F_{jd}/F_d}{F_{(i-1)d}/F_d} \right\rf \\
& = & \left\lf \frac{(F_{id}-F_{(i-1)d})/F_d - 1}{F_{(i-1)d}/F_d} \right\rf \\
& = & \left\lf \frac{F_{id}-F_d}{F_{(i-1)d}} \right\rf - 1 \\
& = & \left\lf \frac{L_d F_{(i-1)d} - F_{(i-2)d} - F_d}{F_{(i-1)d}} \right\rf - 1 \\
& = & L_d-2. 
\end{eqnarray*} 
\vskip 5pt

\noindent We now prove part (iv). In case (i), using eqn.~\eqref{V_identity}, we have  
\begin{eqnarray*}
s(L_m-1) & = & \sum_{i=1}^k {\l}_i V_{n+id} \\
& = & (L_d-2) \sum_{i=1}^{q-1} V_{n+id} + V_{n+(q-1)d} + (L_{d-r}-1) V_{n+qd} + L_r V_{n+(q+1)d} \\
& = & \left( V_{n+qd} - V_{n+(q-1)d} - V_{n+d} + V_n \right) + V_{n+(q-1)d} + (L_{d-r}-1) V_{n+qd} + L_r V_{n+(q+1)d} \\
& = & \left( L_r V_{n+(q+1)d} + L_{d-r} V_{n+qd} \right) - V_{n+d} + V_n \\ 
& = & \left( F_{r+1}+F_{r-1} \right) V_{n+(q+1)d} + \left( F_{d-r+1}+F_{d-r-1} \right) V_{n+qd} - V_{n+d} + V_n \\
& = & \left( F_{r+1} V_{n+(q+1)d} + F_{d-r-1} V_{n+qd} \right) + \left( F_{r-1} V_{n+(q+1)d} + F_{d-r+1} V_{n+qd} \right) - V_{n+d} + V_n \\
& = & F_d \left( V_{n+qd+r+1} + V_{n+qd+r-1} \right) - V_{n+d} + V_n.    
\end{eqnarray*}
In case (ii), using eqn.~\eqref{V_identity}, we have  
\begin{eqnarray*}
s(L_m-1) & = & \sum_{i=1}^k {\l}_i V_{n+id} \\
& = & (L_d-2) \sum_{i=1}^{q-1} V_{n+id} + (L_d-L_{d-r}-1) V_{n+qd} + (L_r-1) V_{n+(q+1)d} \\
& = & \left( V_{n+qd} - V_{n+(q-1)d} - V_{n+d} + V_n \right) + (L_d-L_{d-r}-1) V_{n+qd} + (L_r-1) V_{n+(q+1)d} \\
& = & \left( L_d V_{n+qd} - V_{n+(q-1)d} \right) - V_{n+d} + V_n - L_{d-r} V_{n+qd} +  (L_r-1) V_{n+(q+1)d} \\
& = & V_{n+(q+1)d} - V_{n+d} + V_n - L_{d-r} V_{n+qd} +  (L_r-1) V_{n+(q+1)d} \\
& = & \left( L_r V_{n+(q+1)d} - L_{d-r} V_{n+qd} \right) - V_{n+d} + V_n \\
& = & \left( F_{r+1}+F_{r-1} \right) V_{n+(q+1)d} - \left( F_{d-r+1}+F_{d-r-1} \right) V_{n+qd} - V_{n+d} + V_n \\
& = & \left( F_{r+1} V_{n+(q+1)d} - F_{d-r-1} V_{n+qd} \right) + \left( F_{r-1} V_{n+(q+1)d} - F_{d-r+1} V_{n+qd} \right) - V_{n+d} + V_n \\
& = & F_d \left( V_{n+qd+r+1} + V_{n+qd+r-1} \right) - V_{n+d} + V_n.   
\end{eqnarray*}
In case (iii), using Proposition \ref{F_identity}, part (ii) and eqn.~\eqref{V_identity}, we have  
\begin{eqnarray*}
s(L_m-1) & = & \sum_{i=1}^k {\l}_i V_{n+2i} \\
& = & (L_d-2) \sum_{i=1}^{q-1} V_{n+2i} + (L_d-1) V_{n+qd} - V_{n+(q-1)d} \\
& = & \left( V_{n+qd} - V_{n+(q-1)d} - V_{n+d} + V_n \right) + (L_d-1) V_{n+qd} - V_{n+(q-1)d} \\
& = & L_d V_{n+qd} - L_0 V_{n+(q-1)d} - V_{n+d} + V_n \\
& = & \left( F_{d+1}+F_{d-1} \right) V_{n+qd} - \left( F_1+F_1 \right) V_{n+(q-1)d} - V_{n+d} + V_n \\
& = & \left( F_{d+1} V_{n+qd} - F_1 V_{n+(q-1)d} \right) + \left( F_{d-1} V_{n+qd} - F_1 V_{n+(q-1)d} \right) - V_{n+d} + V_n \\  
& = & F_d \left(V_{n+qd+1} + V_{n+qd-1} \right) - V_{n+d} + V_n.      
\end{eqnarray*}
In each case, $s(L_m-1)=F_d \left(V_{n+m+1}+V_{n+m-1} \right) -V_{n+d}+V_n$. 
\end{Pf}
\vskip 5pt

\subsection{Ap\'ery Set} \label{apery_subsection}
\vskip 5pt

\begin{thm} \label{Greedy_is_optimum}
For any sequence ${\a}_1, \ldots, {\a}_m$ of nonnegative integers, not all zero, 
\[ s\left( \sum_{i=1}^m {\a}_i \frac{F_{id}}{F_d} \right) \le \sum_{i=1}^m {\a}_i V_{n+id}. \]
\end{thm}

\begin{Pf}
We induct on the sum $\s=\sum_{i=1}^m {\a}_i F_{id}/F_d$, where we may assume ${\a}_m \ne 0$ without loss of generality. If $\s=1$, then $m={\a}_1=1$ and the two sides are equal. For some positive integer $\s$, assume the result holds whenever the sum $\sum_{i=1}^m {\a}_i F_{id}/F_d<\s$. 
\vskip 5pt

\noindent Let $F_{kd}/F_d \le \s < F_{(k+1)d}/F_d$ and let ${\l}_1,\ldots,{\l}_k = \text{\sc Greedy}(1,F_{2d}/F_d,F_{3d}/F_d,\ldots,F_{kd}/F_d;\s)$. Suppose ${\a}_1,\ldots,{\a}_m$ is any sequence of nonegative integers such that $\s=\sum_{i=1}^m {\a}_i F_{id}/F_d$; we may assume that ${\a}_m \ge 1$. Note that $m \le k$, for if $m>k$, then $\sum_{i=1}^m {\a}_i F_{id}/F_d \ge F_{(k+1)d}/F_d>\s$. 
\vskip 5pt

\noindent If $m=k$, then $1 \le {\a}_k \le {\l}_k$. By Induction Hypothesis, 
\[ s\left( \sum_{i=1}^k {\a}_i \frac{F_{id}}{F_d} - \frac{F_{kd}}{F_d} \right) \le \sum_{i=1}^k {\a}_i V_{n+id} - V_{n+kd}. \]
Since  
\[ s\left( \sum_{i=1}^k {\a}_i \frac{F_{id}}{F_d} - \frac{F_{kd}}{F_d} \right) =  s\left( \sum_{i=1}^k {\l}_i \frac{F_{id}}{F_d} - \frac{F_{kd}}{F_d} \right) = \sum_{i=1}^k {\l}_i V_{n+id} - V_{n+kd} 
   = s\left( \sum_{i=1}^k {\l}_i \frac{F_{id}}{F_d} \right) - V_{n+kd}, 
\]
we have 
\[ s\left( \sum_{i=1}^k {\a}_i \frac{F_{id}}{F_d} \right) = s\left( \sum_{i=1}^k {\l}_i \frac{F_{id}}{F_d} \right) \le \sum_{i=1}^k {\a}_i V_{n+id}. \]
This proves the Proposition when $m=k$. 
\vskip 5pt

\noindent Suppose $m<k$. By Induction Hypothesis, 
\begin{equation} \label{IH}
s\left( \sum_{i=1}^m {\a}_i \frac{F_{id}}{F_d} - \frac{F_{md}}{F_d} \right) \le \sum_{i=1}^m {\a}_i V_{n+id} - V_{n+md}. 
\end{equation}
Two cases arise: (I) $\sum_{i=1}^m {\a}_i F_{id} - F_{md} \ge F_{kd}$, and (II) $\sum_{i=1}^m {\a}_i F_{id} - F_{md} < F_{kd}$. 
\vskip 5pt

\noindent {\sc Case} (I): Let ${\l}_1^{\pr}, \ldots, {\l}_k^{\pr} = \text{\sc Greedy}\left(1,F_{2d}/F_d,F_{3d}/F_d,\ldots,F_{kd}/F_d;\s-F_{md}/F_d\right)$. Then 
\begin{equation} \label{m<k_subcase1}
s\left( \sum_{i=1}^m {\a}_i \frac{F_{id}}{F_d} - \frac{F_{md}}{F_d} \right) = s\left( \sum_{i=1}^k {\l}_i^{\pr} \frac{F_{id}}{F_d} \right) = \sum_{i=1}^k {\l}_i^{\pr} V_{n+id}. 
\end{equation}
If we replace ${\l}_m^{\pr}$ by ${\l}_m^{\pr}+1$ and retain the other ${\l}_i^{\pr}$, and apply the case $m=k$ discussed above, we get 
\[ s\left( \sum_{i=1}^m {\a}_i \frac{F_{id}}{F_d} \right) =  s\left( \sum_{i=1}^k {\l}_i^{\pr} \frac{F_{id}}{F_d} + \frac{F_{md}}{F_d} \right) \le \sum_{i=1}^k {\l}_i^{\pr} V_{n+id} + V_{n+md} \le 
   \sum_{i=1}^m {\a}_i V_{n+id}. 
\]
from eqn.~\eqref{IH} and eqn.~\eqref{m<k_subcase1}. This proves Case (I). 
\vskip 5pt

\noindent {\sc Case} (II): Since $\s-F_{md}/F_d \ge F_{kd}/F_d-F_{(k-1)d}/F_d>F_{(k-1)d}/F_d$, we have 
\[ {\l}_1^{\pr}, \ldots, {\l}_{k-1}^{\pr} = \text{\sc Greedy}\left(1,F_{2d}/F_d,F_{3d}/F_d,\ldots,F_{(k-1)d}/F_d;\s-F_{md}/F_d\right). \]
\vskip 5pt

\noindent Note that $\s-F_{(k-1)d}/F_d$ lies between $F_{(k-1)d}/F_d$ and $F_{kd}/F_d$. Let 
\[ {\l}_1^{\pr\pr}, \ldots, {\l}_{k-1}^{\pr\pr} = \text{\sc Greedy}\left(1,F_{2d}/F_d,F_{3d}/F_d,\ldots,F_{(k-1)d}/F_d;\s-F_{(k-1)d}/F_d\right). \]
We claim that one of the following cases must arise: (i) ${\l}_i^{\pr\pr}=L_d-2$ for $i \in \{1,\ldots,k-1\}$; (ii) there exists $r \in \{1,\ldots,k-1\}$ such that ${\l}_r^{\pr\pr}=L_d-1$ and ${\l}_i^{\pr\pr}=L_d-2$ for $i \in \{r+1,\ldots,k-1\}$. 
\vskip 5pt

\noindent If neither of these cases is true, then there must exist $t \in \{1,\ldots,k-1\}$ such that ${\l}_t^{\pr\pr} < L_d-2$ and ${\l}_i^{\pr\pr}=L_d-2$ for $i \in \{t+1,\ldots,k-1\}$. But then 
\begin{eqnarray*} 
{\l}_t^{\pr\pr} & = & \left\lf \frac{\s-F_{(k-1)d}/F_d-\sum_{i=t+1}^{k-1} {\l}_i^{\pr\pr} F_{id}/F_d}{F_{td}/F_d} \right\rf \\
& \ge & \left\lf \frac{F_{kd}-F_{(k-1)d}-(L_d-2) \sum_{i=t+1}^{k-1} F_{id}}{F_{td}} \right\rf \\
& \ge & \left\lf \frac{F_{kd}-F_{(k-1)d}-\big(F_{kd}-F_{(k-1)d}-F_{(t+1)d}+F_{td}\big)}{F_{td}} \right\rf \\
& = & \left\lf \frac{F_{(t+1)d}-F_{td}}{F_{td}} \right\rf \\
& \ge & L_d-2
\end{eqnarray*}
by Proposition \ref{F_identity}, parts (iv), (v). This contradiction proves the claim. 
\vskip 5pt

\noindent In Case (i), using Proposition \ref{F_identity}, part (v), we have
\[ \sum_{i=1}^{k-1} {\l}_i^{\pr\pr} \frac{F_{id}}{F_d} = (L_d-2) \sum_{i=1}^{k-1} \frac{F_{id}}{F_d} = \frac{F_{kd}}{F_d} - \frac{F_{(k-1)d}}{F_d} - 1 < \s - \frac{F_{(k-1)d}}{F_d}, \]
contradicting the fact that ${\l}_1^{\pr\pr}, \ldots, {\l}_{k-1}^{\pr\pr}$ is the sequence determined by the Greedy Algorithm for $\s-F_{(k-1)d}$ and $\s \ge F_{kd}/F_d$. This rules out Case (i).  
\vskip 5pt

\noindent In Case (ii), using Proposition \ref{F_identity}, part (vi), we get
\begin{eqnarray} \label{m<k_subcase2}
s\left( \sum_{i=1}^{k-1} {\l}_i^{\pr\pr} \frac{F_{id}}{F_d} \right) + V_{n+(k-1)d} & = & \sum_{i=1}^{k-1} {\l}_i^{\pr\pr} V_{n+id} + V_{n+(k-1)d} \nonumber \\
& = & \sum_{i=1}^{r-1} {\l}_i^{\pr\pr} V_{n+id} + (L_d-2) \sum_{i=r}^{k-1} V_{n+id} + V_{n+rd} + V_{n+(k-1)d} \nonumber \\
& = & \sum_{i=1}^{r-1} {\l}_i^{\pr\pr} V_{n+id} + \big( V_{n+kd} - V_{n+(k-1)d} - V_{n+rd} + V_{n+(r-1)d} \big) \nonumber \\
& & + V_{n+rd} + V_{n+(k-1)d} \nonumber \\
& = & \sum_{i=1}^{r-1} {\l}_i^{\pr\pr} V_{n+id} + V_{n+kd} + V_{n+(r-1)d}. 
\end{eqnarray}
\vskip 5pt 

\noindent Using Proposition \ref{F_identity}, part (v), we have 
\begin{eqnarray*}
\s - \frac{F_{(k-1)d}}{F_d} & = & \sum_{i=1}^{k-1} {\l}_i^{\pr\pr} \frac{F_{id}}{F_d} \\
& = & \sum_{i=1}^{r-1} {\l}_i^{\pr\pr} \frac{F_{id}}{F_d} + (L_d-2) \sum_{i=r}^{k-1} \frac{F_{id}}{F_d} + \frac{F_{rd}}{F_d} \\
& = & \sum_{i=1}^{r-1} {\l}_i^{\pr\pr} \frac{F_{id}}{F_d} + \frac{F_{kd} - F_{(k-1)d} - F_{rd} + F_{(r-1)d}}{F_d} + \frac{F_{rd}}{F_d} \\
& = & \sum_{i=1}^{r-1} {\l}_i^{\pr\pr} \frac{F_{id}}{F_d} + \frac{F_{kd} - F_{(k-1)d} + F_{(r-1)d}}{F_d}. 
\end{eqnarray*}
\vskip 5pt

\noindent By the Induction Hypothesis, 
\[ s\left( \sum_{i=1}^{k-1} {\l}_i^{\pr} \frac{F_{id}}{F_d} - \frac{F_{(k-1)d}}{F_d} + \frac{F_{md}}{F_d} \right) \le \sum_{i=1}^{k-1} {\l}_i^{\pr} V_{n+id} - V_{n+(k-1)d} + V_{n+md}. \]
Applying the case $m=k$ discussed above to $s\left(\s\right)$ and using eqn.~\eqref{m<k_subcase2}, we have  
\begin{eqnarray*} 
s\left( \sum_{i=1}^{r-1} {\l}_i^{\pr\pr} \frac{F_{id}}{F_d} + \frac{F_{kd}}{F_d} + \frac{F_{(r-1)d}}{F_d} \right) & \le & \sum_{i=1}^{r-1} {\l}_i^{\pr\pr} V_{n+id} + V_{n+kd} + V_{n+(r-1)d} \\
& = & s\left( \sum_{i=1}^{k-1} {\l}_i^{\pr\pr} \frac{F_{id}}{F_d} \right) + V_{n+(k-1)d} \\
& = & s\left( \sum_{i=1}^{k-1} {\l}_i^{\pr} \frac{F_{id}}{F_d} - \frac{F_{(k-1)d}}{F_d} + \frac{F_{md}}{F_d} \right) + V_{n+(k-1)d} \\
& \le & \sum_{i=1}^{k-1} {\l}_i^{\pr} V_{n+id} + V_{n+md} \\
& = & s\left( \sum_{i=1}^m {\a}_i \frac{F_{id}}{F_d} - \frac{F_{md}}{F_d} \right) + V_{n+md} \\
& \le & \sum_{i=1}^m {\a}_i V_{n+id}.   
\end{eqnarray*} 

\noindent This completes Case (ii), and the proof. 
\end{Pf}
\vskip 5pt

\begin{lem} \label{s_increasing}
For any positive integer $m$, $s(m)<s(m+1)$. 
\end{lem}

\begin{Pf}
We induct on $m$. Note that $V_{n+d}=s(1)<s(2)=2V_{n+d}$. Assume $s(i-1)<s(i)$ for $1 \le i \le m$. If $m=F_{kd}/F_d-1$ for some $k$, then $s(m)=V_{n+kd}-V_{n+d}+V_n<V_{n+kd}=s(m+1)$ by Proposition \ref{greedy_algo_for_F_kd}, part (iii). Otherwise $F_{kd}/F_d \le m<F_{(k+1)d}/F_d-1$, and so $s(m)=s(m-F_{kd}/F_d)+V_{n+kd}$ while $s(m+1)=s(m+1-F_{kd}/F_d)+V_{n+kd}$. By Induction Hypothesis, $s(m-F_{kd}/F_d)<s(m+1-F_{kd}/F_d)$, so that $s(m)<s(m+1)$, proving the Proposition by induction. 
\end{Pf}
\vskip 5pt

\begin{thm} \label{Apery_set}
Let $\gcd(V_1,V_2)=\gcd(V_n,F_d)=1$, where $d$ is even. The Ap\'{e}ry set for $S=\la V_n, V_{n+d}, V_{n+2d}, \ldots \ra$ is given by  
\[ \Ap(S,V_n) = \{ s(x): 1 \le x \le V_n-1 \} \cup \{0\}. \] 
\end{thm}

\begin{Pf}
For $x \in \{1,\ldots,V_n-1\}$, we show that $s(x)$ is the least positive integer in $S$ that is congruent to $V_{n+d}\,x$ modulo $V_n$. This proves the result since $\{V_{n+d}\,x: 1 \le x \le V_n-1\}$ is the set of non-zero residues modulo $V_n$ as $\gcd(V_n,V_{n+d})=1$. 

Suppose $s \in S$ is congruent to $V_{n+d}\,x$ modulo $V_n$. Then $s=\sum_{i \ge 0} {\a}_i V_{n+id}$, with each ${\a}_i \ge 0$. Since $s \equiv V_{n+d}\,x\pmod{V_n}$ and $V_{n+id} \equiv V_{n+d} F_{id}/F_d\pmod{V_n}$ by Proposition \ref{F_identity}, part (ii), we have $\sum_{i \ge 1} {\a}_i F_{id}/F_d \equiv x\pmod{V_n}$ as $\gcd(V_n,V_{n+d})=1$. Since $x \le V_n-1$, we have $x \le \sum_{i \ge 1} {\a}_i F_{id}/F_d$, so that  
\[ s(x) \le s\left( \sum_{i \ge 1} {\a}_i \frac{F_{id}}{F_d} \right) \le \sum_{i \ge 1} {\a}_i V_{n+id} \le s \]
by Theorem \ref{Greedy_is_optimum} and Lemma \ref{s_increasing}. 
\end{Pf}
\vskip 5pt

\subsection{The Frobenius number and Genus in Some Special Cases} \label{spcases_subsection}
\vskip 5pt

\begin{thm} \label{Frob_gen}
Let $\gcd(V_1,V_2)=\gcd(V_n,F_d)=1$, where $d$ is even. If $S=\la V_n, V_{n+d}, V_{n+2d}, \ldots \ra$, then 
\begin{itemize}
\item[{\rm (i)}]
\[ \F(S) = s(V_n-1) - V_n, \]
\item[{\rm (ii)}]
\[ \g(S) = \frac{1}{V_n} \left( \sum_{x=1}^{V_n-1} s(x) \right) - \frac{V_n-1}{2}. \] 
\end{itemize}
\end{thm}

\begin{Pf}
These are direct consequences of Proposition \ref{prelims}, Theorem \ref{Apery_set} and Lemma \ref{s_increasing}. 
\end{Pf}
\vskip 5pt

\begin{cor} \label{Frob_special}
\begin{itemize}
\item[]
\item[{\rm (i)}]
If $S_1=\la F_n, F_{n+d}, F_{n+2d}, \ldots \ra$, $n \ge 3$, then $\F(S_1)=F_d F_{2n}-F_{n+d}$. 
\item[{\rm (ii)}]
If $S_2=\la L_n, L_{n+d}, L_{n+2d}, \ldots \ra$, $n \ge 4$, then $\F(S_2)=F_d \left( L_{2n+1}+L_{2n-1} \right)-L_{n+d}$. 
\end{itemize}
\end{cor}

\begin{Pf}
This is a direct consequence of Theorem \ref{Frob_gen} and Propositions \ref{greedy_fibonacci}, \ref{greedy_lucas}. 
\end{Pf}
\vskip 5pt

\noindent Computation of $\g(S)$ is difficult in the general case. In the following result we compute the genus in the special case of Fibonacci and Lucas subsequences. The result is in terms of the $k^{\text{th}}$ term of sequences that jointly use first order recurrences, and that can be explicitly solved. 
\vskip 5pt

\begin{prop} \label{genus_recurrence}
For $k \ge 1$, define $A_k$ and $B_k$ as follows: 
\begin{equation}
A_k = \sum_{x=1}^{\tfrac{F_{(k+1)d}}{F_d}-1}\: s(x), \quad B_k = \sum_{x=1}^{\tfrac{F_{(k+1)d}-F_{kd}}{F_d}-1}\: s(x).  
\end{equation}
Then $A_k$ and $B_k$ satisfy the joint first order recurrences given by  
\begin{eqnarray}
A_{k+1} & = & (L_d-1) A_k + B_k + (L_d-1) V_{n+(k+1)d} \frac{F_{(k+2)d} - F_{kd}}{2F_d}; \\[5pt]
B_{k+1} & = & (L_d-2) A_k + B_k + (L_d-2) V_{n+(k+1)d} \frac{F_{(k+2)d} - F_{(k+1)d - F_{kd}}}{2F_d}. 
\end{eqnarray}
with $A_1=\frac{1}{2}V_{n+d} L_d(L_d-1)$ and $B_1=\frac{1}{2}V_{n+d} (L_d-1)(L_d-2)$. 
\end{prop}

\begin{Pf}
Let $x$ be a positive integer. Then $F_{kd}/F_d \le x < F_{(k+1)d}/F_d$ for some positive integer $k$. Thus,   
\[ i \frac{F_{kd}}{F_d} \le x < (i+1) \frac{F_{kd}}{F_d} \]
for some $i \in \{1,\ldots,L_d-1\}$. From the definition of $s$, 
\begin{equation} \label{s_reduction}
s(x) = i\,V_{n+kd} + s\left(x-i\,\frac{F_{kd}}{F_d}\right). 
\end{equation}
\vskip 5pt

\noindent Therefore, from eqn.~\eqref{s_reduction} we have  
\begin{eqnarray}
\sum_{x=1}^{i\,\frac{F_{kd}}{F_d}-1} s(x) & = & \sum_{x=1}^{\frac{F_{kd}}{F_d}-1} s(x) + \sum_{x=\frac{F_{kd}}{F_d}}^{2\,\frac{F_{kd}}{F_d}-1} s(x) + \cdots + \sum_{x=\frac{(i-1)\,F_{kd}}{F_d}}^{i\,\frac{F_{kd}}{F_d}-1} s(x) \nonumber \\
& = & A_{k-1} + \left( \tfrac{F_{kd}}{F_d} V_{n+kd} + A_{k-1} \right) + \left( \tfrac{2F_{kd}}{F_d} V_{n+kd} + A_{k-1} \right) + \cdots + \left( \tfrac{(i-1)F_{kd}}{F_d} V_{n+kd} + A_{k-1} \right) \nonumber \\[5pt]
& = & i\,A_{k-1} + \frac{i(i-1)}{2} \frac{F_{kd}}{F_d} V_{n+kd}. \label{i_sum}
\end{eqnarray}
\vskip 5pt

\noindent We have 
\begin{equation} \label{fibluc_identity}
F_{(k+2)d} = (L_d-1) F_{(k+1)d} + F_{(k+1)d} - F_{kd} 
\end{equation}
by Proposition \ref{F_identity}, part (iii). 
\vskip 5pt

\noindent Hence, from eqn.~\eqref{s_reduction} and eqn.~\eqref{i_sum}, and using eqn.~\eqref{fibluc_identity}, we have  
\begin{eqnarray*}
A_{k+1} & = & \sum_{x=1}^{\frac{F_{(k+2)d}}{F_d}-1} s(x) \\
& = & \sum_{x=1}^{(L_d-1)\frac{F_{(k+1)d}}{F_d}-1} s(x) + \sum_{x=(L_d-1)\frac{F_{(k+1)d}}{F_d}}^{\frac{F_{(k+2)d}}{F_d}-1} s(x) \\
& = & \left( (L_d-1) A_k + \tfrac{(L_d-2)(L_d-1)}{2} \tfrac{F_{(k+1)d}}{F_d} V_{n+(k+1)d} \right) + (L_d-1) \left( \tfrac{F_{(k+1)d}}{F_d} - \tfrac{F_{kd}}{F_d} \right) V_{n+(k+1)d} \\
& & + \sum_{x=0}^{\frac{F_{(k+1)d}-F_{kd}}{F_d}-1} s(x) \\
& = & (L_d-1) A_k + B_k + \frac{L_d-1}{2F_d} V_{n+(k+1)d} \left( L_d\,F_{(k+1)d} - 2F_{kd} \right) \\[5pt] 
& = & (L_d-1) A_k + B_k + (L_d-1) V_{n+(k+1)d} \frac{F_{(k+2)d} - F_{kd}}{2F_d}. 
\end{eqnarray*}
\vskip 5pt

\noindent The derivation of the formula for $B_{k+1}$ follows along similar lines. We have 
\begin{eqnarray*}
B_{k+1} & = & \sum_{x=1}^{\frac{F_{(k+2)d}-F_{(k+1)d}}{F_d}-1} s(x) \\
& = & \sum_{x=1}^{(L_d-2)\frac{F_{(k+1)d}}{F_d}-1} s(x) + \sum_{x=(L_d-2)\frac{F_{(k+1)d}}{F_d}}^{\frac{F_{(k+2)d}-F_{(k+1)d}}{F_d}-1} s(x) \\
& = & \left( (L_d-2) A_k + \tfrac{(L_d-3)(L_d-2)}{2} \tfrac{F_{(k+1)d}}{F_d} V_{n+(k+1)d} \right) + (L_d-2) \left( \tfrac{F_{(k+1)d}}{F_d} - \tfrac{F_{kd}}{F_d} \right) V_{n+(k+1)d} \\
& & + \sum_{x=0}^{\frac{F_{(k+1)d}-F_{kd}}{F_d}-1} s(x) \\
& = & (L_d-2) A_k + B_k + \frac{L_d-2}{2F_d} V_{n+(k+1)d} \left( (L_d-1)F_{(k+1)d} - 2F_{kd} \right) \\[5pt] 
& = & (L_d-2) A_k + B_k + (L_d-2) V_{n+(k+1)d} \frac{F_{(k+2)d} - F_{(k+1)d} - F_{kd}}{2F_d}. 
\end{eqnarray*}
\vskip 5pt

\noindent Moreover, 
\[ A_1 = \sum_{x=1}^{\frac{F_{2d}}{F_d}-1} s(x) = \sum_{x=1}^{L_d-1} V_{n+d}\,x = \frac{1}{2}(L_d-1)L_d V_{n+d} \] 
and 
\[ B_1 = \sum_{x=1}^{\frac{F_{2d}-F_d}{F_d}-1} s(x) = \sum_{x=1}^{L_d-2} V_{n+d}\,x = \frac{1}{2}(L_d-2)(L_d-1) V_{n+d}. \] 
\end{Pf}
\vskip 5pt

\noindent Let $X$ be a positive integer and let $k$ be such that $F_{kd}/F_d \le X < F_{(k+1)d}/F_d$. Let 
\[ {\l}_1,\ldots,{\l}_k = \text{\sc Greedy}(1,F_{2d}/F_d,F_{3d}/F_d,\ldots,F_{kd}/F_d;X). \]
Assume 
\[ {\l}_i = \begin{cases} 
                 L_d-2 & \mbox{ if } 1 \le i \le k-2; \\
                 b & \mbox{ if } i=k-1; \\
                 a & \mbox{ if } i=k, 
               \end{cases}
\]
where $a, b \le L_{d-1}$ with $(a,b) \ne (L_{d-1},L_{d-1})$. Then 
\begin{eqnarray}
\sum_{x=1}^X s(x) & = & \sum_{x=1}^{a\,\frac{F_{kd}}{F_d}-1} s(x) + \sum_{x=a\,\frac{F_{kd}}{F_d}}^{\frac{a\,F_{kd}+b\,F_{(k-1)d}}{F_d}-1} s(x) + \sum_{x=\frac{a\,F_{kd}+b\,F_{(k-1)d}}{F_d}}^X s(x) \nonumber \\
& = & a A_{k-1} + \tfrac{(a-1)a}{2} \tfrac{F_{kd}}{F_d} V_{n+kd} + ab \tfrac{F_{(k-1)d}}{F_d} V_{n+kd} + b A_{k-2} + \tfrac{(b-1)b}{2} \tfrac{F_{(k-1)d}}{F_d} V_{n+(k-1)d} \nonumber \\
& & + \tfrac{F_{(k-1)d}-F_{(k-2)d}}{F_d} \left( a V_{n+kd} + b V_{n+(k-1)d} \right) + B_{k-2} \nonumber \\
& = & a A_{k-1} + B_{k-2} + b A_{k-2} + \tfrac{a}{F_d} \left( \tfrac{a-1}{2} F_{kd} + (b+1) F_{(k-1)d} - F_{(k-2)d} \right) V_{n+kd} \nonumber \\
& & + \tfrac{b}{F_d} \left( \tfrac{b+1}{2} F_{(k-1)d} - F_{(k-2)d} \right) V_{n+(k-1)d}. \label{ugly_sum} 
\end{eqnarray}
\vskip 5pt

\noindent Recall from Theorem \ref{Frob_gen}, part (ii) that $\g(S)$ may be determined from $\sum_{x=1}^{V_n-1} s(x)$. When $V_n=F_n$ or $L_n$, the ${\l}_i$'s for $F_n-1$ and for $L_n-1$ are of the form given in the above discussion in one of the cases; in the other cases, there is the presence of an additional constant corresponding to ${\l}_{k-2}$ which is distinct from $L_d-2$. In such cases, as the ${\l}_i$'s take the above form, eqn.~\eqref{ugly_sum} provides a closed form expression for $\sum_{x=1}^{V_n-1} s(x)$. A similar expression may also be derived in case ${\l}_k, {\l}_{k-1}, {\l}_{k-2}$ are all distinct from $L_d-2$. We remark that the expression derived in eqn.~\eqref{ugly_sum} involves the terms from the sequences $A_k$ and $B_k$. In principle, these may be evaluated by solving the two recurrences in Proposition \ref{genus_recurrence}.


\vskip 20pt

\begin{thebibliography}{99}

\bibitem{BKT15}
S. S. Batra, N. Kumar and A. Tripathi, On a linear Diophantine problem involving the Fibonacci and Lucas sequences, {\it Integers\/} {\bf 15} (2015), Article A26, 12 pp. 

\bibitem{BS62}
A. Brauer and J. E. Shockley, On a problem of Frobenius, {\it J. Reine Angew. Math.\/} {\bf 211} (1962), 215-220.
 
\bibitem{Cur90}
F. Curtis, On formulas for the Frobenius number of a numerical semigroup, {\it Math. Scand.\/} {\bf 67} (1990), 190-92. 

\bibitem{MRR07}
J. M. Mar\'{i}n, J. L. Ram\'{i}rez Alfons\'{ı}n and M. P. Revuelta, On the Frobenius number of Fibonacci Numerical Semigroups, {\it Integers\/} {\bf 7} (2007), Article A14, 7 pp. 

\bibitem{Mat09}
G. L. Matthews, Frobenius Numbers of Generalized Fibonacci Semigroups, {\it Integers\/} {\bf 9\,Supplement} (2009), Article 9, 7 pp. 

\bibitem{PRT24}
S. Panda, K. Rai and A. Tripathi, On the Frobenius Problem for Some Generalized Fibonacci Subsequences - I, preprint. 

\bibitem{Ram05}
J. L. Ram\'{i}rez Alfons\'{ı}n, The Diophantine Frobenius Problem, {\it Oxford University Press\/}, 2005, 259 pp. 

\bibitem{RG-S09}
J. C. Rosales and P. A. Garc\'{i}a-S\'{a}nchez, Numerical Semigroups, {\it Springer-Verlag\/}, 2009, 181 pp. 

\bibitem{Ryb24}
D. Rybin, When greedy gives optimal: A unified approach, {\it Discrete Optim.\/} {\bf 51} (2024), 100824. 

\bibitem{Sel77}
E. S. Selmer, On the linear diophantine problem of Frobenius, {\it J. Reine Angew. Math.\/} {\bf 293/294} (1977), 1-17.



\end{thebibliography}
\end{document}